\documentclass{elsarticle-modified}
\usepackage{amssymb,amsmath,amsfonts}
\usepackage{bm}
\usepackage[utf8]{inputenc}
\usepackage[T1]{fontenc}
\usepackage{xcolor}
\usepackage{mathrsfs}
\usepackage{cancel}
\usepackage{nicefrac}
\usepackage{boxedminipage}

\numberwithin{equation}{section}

\definecolor{light-gray}{gray}{0.95}

\journal{J. Comput. Appl. Math.}

\begin{document}

\newcommand\rev[1]{{#1}}
\newenvironment{remark}{{\bf Remark.}\/}{}
\newtheorem{theorem}{Theorem}
\newtheorem{proposition}[theorem]{Proposition}
\newtheorem{lemma}[theorem]{Lemma}

\newcommand\sig{\sigma}

\newcommand{\dd}{{\mathrm d}}
\newcommand\dt{\tfrac{\dd}{\dd t}}
\newcommand\dtau{\tfrac{\dd}{\dd \tau}}
\newcommand\dtnot{\tfrac{\dd}{\dd t_0}}
\newcommand\du{\tfrac{\dd}{\dd u}}
\newcommand\dunot{\tfrac{\dd}{\dd u_0}}
\newcommand\ds{\tfrac{\dd}{\dd s}}
\newcommand\pds{\tfrac{\partial}{\partial s}}
\newcommand\pdsig{\tfrac{\partial}{\partial \sig}}
\newcommand\pdone{\partial_1}
\newcommand\pdtwo{\partial_2}
\newcommand\pdthree{\partial_3}
\newcommand\pdtau{\tfrac{\partial}{\partial \tau}}
\newcommand\pdt{\tfrac{\partial}{\partial t}}
\newcommand\pdtnot{\tfrac{\partial}{\partial t_0}}
\newcommand\pdu{\tfrac{\partial}{\partial u}}
\newcommand\pdunot{\tfrac{\partial}{\partial u_0}}

\newcommand\ol[1]{\overline{#1}}
\newcommand\nool[1]{#1}
\newcommand\ul[1]{\underline{#1}}
\newcommand\noul[1]{#1}
\newcommand\ub[2]{\underbrace{#1}_{#2}}
\newcommand\noub[2]{#1}

\newcommand{\CCC}{{{\mathbb{C}}\vphantom{|}}}
\newcommand{\NN}{{{\mathbb{N}}\vphantom{|}}}
\newcommand{\RR}{{{\mathbb{R}}\vphantom{|}}}

\renewcommand{\dd}{{\mathrm d}}
\newcommand{\ee}{{\mathrm e}}
\newcommand{\ii}{{\mathrm i}}

\newcommand\clas{{\bm c}}
\newcommand\symm{{\bm s}}
\newcommand\Id{\text{Id}}
\newcommand\Order{{\mathscr{O}}}
\renewcommand\th{\tfrac{1}{2}}
\newcommand\tth{\frac{t}{2}}
\newcommand\nEE[1]{\ee^{\,{#1}}}
\newcommand\nE{{\mathcal E}}
\newcommand\nD{{\mathcal D}}
\newcommand\nDc{{\mathcal D}_{\clas}}
\newcommand\nDs{{\mathcal D}_{\symm}}
\newcommand\nDstilde{{\widetilde{\mathcal D}}_{\symm}}
\newcommand\nL{{\mathcal L}}
\newcommand\nLLtildes{{\widetilde\nL}_\symm}
\newcommand\nS{{\mathcal S}}
\newcommand\nSD{{\widehat\nS}_\symm}
\newcommand\nLD{{\widehat\nL}_\symm}
\newcommand\nP{\widetilde{\mathcal{L}}}
\newcommand\nR{{\mathcal R}}
\newcommand\nX{{\mathcal X}}
\newcommand\Bcheck{\check{B}}

\newcommand\canc[1]{\fcolorbox{gray}{light-gray}{\text{$#1$}}}

\begin{frontmatter}

\title{Symmetrized local error estimators for time-reversible one-step methods
       in nonlinear evolution equations}

\author{Winfried Auzinger}
\ead{w.auzinger@tuwien.ac.at}
\address{Technische Universit{\"a}t Wien,
Institut f{\"u}r Analysis und Scientific Computing,
Wiedner Hauptstrasse 8--10/E101, A-1040 Wien, Austria}
\ead[url]{www.asc.tuwien.ac.at/\~{}winfried}

\author{Harald Hofst{\"a}tter}
\address{Universit{\"a}t Wien,
Institut f{\"u}r Mathematik,
Oskar-Morgenstern-Platz 1, A-1090 Wien, Austria}
\ead{hofi@harald-hofstaetter.at}
\ead[url]{www.harald-hofstaetter.at}

\author{Othmar Koch\corref{CA}}
\address{Universit{\"a}t Wien,
Institut f{\"u}r Mathematik,
Oskar-Morgenstern-Platz 1, A-1090 Wien, Austria}
\ead{othmar@othmar-koch.org}
\ead[url]{www.othmar-koch.org}

\cortext[CA]{Corresponding author.}

\begin{abstract}
Prior work on computable defect-based local error estimators
for (linear) time-reversible integrators is extended to nonlinear and nonautonomous
evolution equations. We prove that the asymptotic results from the linear case
[W.\;Auzinger and O.\;Koch,
An improved local error estimator for symmetric time-stepping schemes,
Appl.\;Math.\;Lett.~82\,(2018), pp.\;106--110]
remain valid, i.e., the modified estimators yield an
improved asymptotic order as the step size goes to zero.
Typically, the computational effort
is only slightly higher than for conventional defect-based
estimators, and it may even be lower in some cases.
We illustrate this by some examples and present numerical
results for evolution equations of Schr{\"o}dinger type,
solved by either time-splitting or Magnus-type integrators.
\rev{Finally, we demonstrate that adaptive time-stepping
schemes can be successfully based on our local error estimators.}
\end{abstract}

\begin{keyword}
Nonlinear evolution equations \sep
numerical time integration \sep
one-step methods \sep
time-reversible schemes \sep
splitting methods \sep
commutator-free Magnus-type methods \sep
Magnus integrators \sep
local error estimation
\MSC[2010] 65L05 \sep 65L20 \sep 65M12
\end{keyword}

\end{frontmatter}


\section{Introduction} \label{sec:intro}

We consider the extension of a defect-based estimator for the local error
of self-adjoint time-stepping schemes of even order $p$, which was introduced in~\cite{auzingeretal18a}
for the linear time-independent case, to nonlinear evolution equations
(we set $ t_0=0 $),
\begin{equation} \label{y'=F(y)}
\dt u(t) = F(u(t)), \quad u(0)=u_0.
\end{equation}
We define a symmetrized version of the defect to serve as the basis for the construction
of a local error estimator in the nonlinear case, thus representing an extension of~\cite{auzingeretal18a}. The error
estimator is derived from a representation
of the local error in terms of the symmetrized defect,
based on a modified nonlinear variation-of-constant formula.
Its deviation from the exact local error is one order in the step-size more precise than an analogous error estimator
based on the classical defect, for the latter see for instance~\cite{auzingeretal12a,auzingeretal13a,auzingeretal13b,auzingeretal14a}.
Our theoretical analysis is based on the assumption that the problem is smooth
(the right-hand side is bounded and differentiable with bounded derivatives as required in
the analysis) with a unique, smooth solution. In this sense, our treatment is formal
and in practical applications with unbounded right-hand side, different techniques
are required to deduce the required regularity assumptions in order to establish
high-order convergence, see for instance~\cite{auzingeretal13b}.

We also point out that in addition to the practical merit of providing a
more precise estimator enabling a better choice of adaptive time-steps
and a higher-order corrected solution if desired, the approach
has potential advantages for theoretical purposes.

In the analysis of local errors and error estimators for self-adjoint schemes,
the representation of the local error in terms
of the symmetrized defect can be rewritten in a way such that
its analysis can be based on an asymptotic expansion
in even powers of the stepsize.
Applications of this type of analysis will be reported elsewhere.

\paragraph{Outline}
In Section~\ref{sec:def} we introduce the notions `classical defect'
and the new `symmetrized defect' associated with
one-step integrators for nonlinear evolution equations
in the autonomous\footnote{The extension to nonautonomous
                           problems is deferred to Section~\ref{sec:nonaut}.}
form~\eqref{y'=F(y)}.
A well-known integral representation of the local error in terms of the classical defect
is obtained from the nonlinear variation-of-constant
formula (V.O.C., also referred to as Gr{\"o}bner-Alexeev-Lemma~\cite{haireretal87}),
this is recapitulated in~Theorem~\ref{th:GAL}. Then,
in Theorem~\ref{th:GALS} we present a modified nonlinear V.O.C.~formula
leading to an integral representation of the local error
in terms of the symmetrized defect.

An Hermite-type quadrature approximation to the ensuing
integral representation provides a computable
defect-based local error estimator, see Section~\ref{sec:locerrest}.
In particular, Theorem~\ref{th:Sdach} shows that the
symmetrized error estimator is asymptotically correct,
and for the case of a self-adjoint scheme
it is of an improved asymptotic quality compared with the analogous
classical estimator.
Here the required regularity of the problem data and of the exact solution is tacitly assumed.

In Sections~\ref{sec:appl-auto} and~\ref{sec:nonaut} we study
the application of these ideas to particular examples of self-adjoint schemes.
In Section~\ref{subsec:imr}, the results are particularized to the implicit midpoint rule
to show a concrete example of an implicit one-step method.
In Section~\ref{subsec:strang},
Strang splitting is discussed, and the algorithmic realization for general splitting
methods is given in Section~\ref{subsec:algsplit}.

In Section~\ref{sec:nonaut}, the nonautonomous case is considered.
In order to illustrate the extension of our ideas to this case,
we give details for linear problems with a $ t $\,-\,dependent
right-hand side.
Section~\ref{subsec:emr} shows the realization for the exponential midpoint rule, and
Section~\ref{subsec:cfm} contains the algorithmic implementation for general commutator-free
Magnus-type and classical Magnus methods.

In Section~\ref{sec:num}, numerical examples for a splitting
approximation to a cubic nonlinear Schr{\"o}dinger equation and Magnus-type exponential
integrators applied to a time-dependent Rosen--Zener model \rev{support} the theoretical results,
\rev{and adaptive time-stepping based on the new error estimator is illustrated}.


\paragraph{Notation and preliminaries}
The flow associated with~\eqref{y'=F(y)}
is denoted by $ \nE(t,u) $, such that the solution of~\eqref{y'=F(y)}
is $ u(t) = \nE(t,u_0) $.
By $ \pdone\nE(t,u_0) $ and $ \pdtwo\nE(t,u_0) $  we denote
the derivatives of $ \nE $ with respect to its first and second
arguments, respectively.
By definition, $ \nE(t,u_0) $ satisfies
\begin{equation*}
\pdone\nE(t,u_0) = F(\nE(t,u_0)), \quad \nE(0,u_0) = u_0.
\end{equation*}
We will repeatedly make use of the following fundamental
identity.\footnote{For the nonautonomous case see Lemma~\ref{lem:fin}
                   in Section~\ref{sec:nonaut}.}

\rev{
\begin{lemma} \label{lem:fi}
\begin{equation} \label{fi}
[\,\pdone\nE(t,u_0) =\,]~
F(\nE(t,u_0)) = \pdtwo\nE(t,u_0) \cdot F(u_0).
\end{equation}
\end{lemma}
{\em Proof.}
\eqref{fi} is a consequence of the first-order variational equation
for $ \nE(t,u) $,
see~\cite[Theorem I.14.3]{haireretal87},~\cite[Appendix~A]{auzingeretal13b}.
The simple direct proof given in~\cite[(3.7)]{DescombesThalhammerLie}
proceeds from the identity
\begin{equation*}
\nE(t+s,u_0) = \nE(t,\nE(s,u_0)).
\end{equation*}
Differentiation with respect to~$ s $ gives
\begin{align*}
\pds\nE(t+s,u_0)
&= \pdone\nE(t+s,u_0), \\
\pds\nE(t+s,u_0)\big|_{s=0}
&= \pdone\nE(t,u_0) = F(\nE(t,u_0)),
\end{align*}
and on the other hand,
\begin{align*}
\pds\nE(t,\nE(s,u_0))
&= \pdtwo\nE(t,\nE(s,u_0)) \cdot \pdone\nE(s,u_0), \\
\pds\nE(t,\nE(s,u_0))\big|_{s=0}
&= \pdtwo\nE(t,u_0) \cdot \pdone\nE(0,u_0) =  \pdtwo\nE(t,u_0) \cdot F(u_0), \\
\end{align*}
which completes the proof.
\hfill $ \square $
}

\section{Classical and symmetrized defects for one-step integrators}
\label{sec:def}
Consider an approximation to the given problem~(\ref{y'=F(y)})
defined by the flow
\begin{equation} \label{S(t_0,u_0)}
\nS(t,u_0) \approx \nE(t,u_0), \quad \nS(0,u_0) = u_0,
\end{equation}
of a consistent one-step scheme with stepsize $ t $,
starting at $ (0,u_0) $.
We assume that the scheme has order~$ p $,
i.e., the local error
\begin{equation} \label{L(t_0,u_0)}
\nL(t,u_0) = \nS(t,u_0) - \nE(t,u_0)
\end{equation}
satisfies $ \nL(t,u_0) = \Order(t^{p+1}) $.

We call
\begin{equation} \label{cdnonlin}
\nDc(t,u) = \pdone\nS(t,u) - F(\nS(t,u)) = \Order(t^p)
\end{equation}
the {\em classical defect}\, associated with $ \nS(t,u) $.
The local error can be represented in terms of the classical defect via
the well-known nonlinear variation-of-constant formula,
the so-called Gr{\"o}bner-Alekseev Lemma.
For convenience we restate this in a form required in our context
and also include the 
proof following\footnote{See~\cite[Figure~I.14.1]{haireretal87},
 {\em Lady Windermere's Fan, Act 2.}}~\cite[Theorem~I.14.5]{haireretal87}
(see also~\cite[Theorem~3.3]{DescombesThalhammerLie}).
We formulate it in a concise way making direct use of~\eqref{fi}.

\begin{theorem}\label{th:GAL}
In terms of the classical defect~\eqref{cdnonlin},
the local error satisfies the integral representation
\begin{equation} \label{GAL-identity}
\nL(t,u_0) =
\int_0^t \pdtwo\nE(t-s,\nS(s,u_0)) \cdot \nDc(s,u_0)\,\dd s.
\end{equation}
\end{theorem}
{\em Proof.}
For fixed $ t $, let
\begin{align*}
y(s) &= \nS(s,u_0), \\
z(s) &= \nE(t-s,y(s)).
\end{align*}
In this notation, we have
\begin{align*}
z(s) &= \nE(t-s,\nS(s,u_0)), \\
\text{satisfying} \quad z(0) &= \nE(t,u_0), ~~z(t) = \nS(t,u_0).
\end{align*}
Thus,
\begin{equation} \label{Lint-1}
\nL(t,u_0) = \nS(t,u_0) - \nE(t,u_0) = \int_0^t \ds z(s) \,\dd s,
\end{equation}
with
\begin{equation*}
\ds z(s) = -F(z(s)) + \pdtwo\nE(t-s,y(s)) \cdot \ds y(s).
\end{equation*}
Now, using~\eqref{fi}\footnote{{\em Mutatis mutandis:} $ s,t-s $ and $ y(s) $
                               play the role of $ 0,t $ and $ u_0 $
                               from~\eqref{fi}.}
this can be rewritten in the form
\begin{align*}
\ds z(s)
&= \ub{-F(\nE(t-s,y(s))) + \pdtwo\nE(t-s,y(s)) \cdot F(y(s))}{=\,0} \\
& \quad {} + \pdtwo\nE(t-s,y(s)) \cdot
  \big( \ds y(s) - F(y(s)) \big) \\
&= \pdtwo\nE(t-s,y(s)) \cdot \nDc(s,u_0),
\end{align*}
and together with~\eqref{Lint-1},
identity~\eqref{GAL-identity} immediately follows.
\hfill $ \square $

\begin{remark}
Due to~\eqref{fi}, an alternative, plausible way to define the defect is
\begin{equation} \label{cdnonlin-alt}
\nD(t,u) = \pdone\nS(t,u) - \pdtwo\nS(t,u) \cdot F(u).
\end{equation}
Then, 
\begin{equation*}
\nL(t,u_0) =
\int_0^t \ds \nS(s,\nE(t-s,u_0))\,\dd s =
\int_0^t \nD(s,\nE(t-s,u_0))\,\dd s.
\end{equation*}
\end{remark}

\begin{remark} \label{rem:DcD}
We can express the modified defect~\eqref{cdnonlin-alt}
in terms of $ \nDc(t,u) $ plus a higher-order perturbation,
\begin{align*}
&\pdone\nS(t,u) - \pdtwo\nS(t,u) \cdot F(u) \\
&= \big( \pdone\nS(t,u) - F(\nS(t,u)) \big)
   + \big( F(\nS(t,u) - \pdtwo\nS(t,u) \cdot F(u) \big) \\
&= \nDc(t,u)
   + \ub{\big( F(\nE(t,u)) - \pdtwo\nE(t,u) \cdot F(u) \big)}{=\,0} \\
& \quad {} + \ub{\big( F(\nS(t,u)) - F(\nE(t,u)) \big)}
                {=\,\Order(t^{p+1})}
 + \ub{\big( \pdtwo\nS(t,u) - \pdtwo\nE(t,u) \big)}
                    {=\,\Order(t^{p+1})} \cdot\,F(u) \\
&= \nDc(t,u) + \Order(t^{p+1}).
\end{align*}
\end{remark}

Also, e.g., a convex combination of~\eqref{cdnonlin}
and~\eqref{cdnonlin-alt} represents a plausible defect.
In particular, we will consider the arithmetic mean of~\eqref{cdnonlin}
and~\eqref{cdnonlin-alt} (see~\eqref{sdnonlin} below),
and we will introduce a symmetrized variant of Theorem~\ref{th:GAL},
see Theorem~\ref{th:GALS} below.

\subsection{Symmetrization}
The following considerations are relevant for the case
where the approximate flow $ \nS $ is self-adjoint
(symmetric, time-reversible),\footnote{Definition~\eqref{sdnonlin} and
                            the assertion of Theorem~\ref{th:GALS}
                            are independent of this assumption.
                            However, our results derived later on
                            essentially depend on it, in particular
                            Theorem~\ref{th:Sdach}.}
i.e.,
\begin{equation} \label{Ssymm}
\nS(-t,\nS(t,u)) = u.
\end{equation}
Self-adjoint schemes have an even order $ p $,
see~\cite[Theorem~II.3.2]{haireretal02}.

The identity\footnote{
  In the terminology of Lie calculus
  (cf.~for instance~\cite{haireretal02}), with
  \begin{equation*}
  (D_F\,G)(u) := G'(u) \cdot F(u) =
  \tfrac{\dd}{\dd t} G(\nE(t,u))|_{t=0},
  \end{equation*}
  and
  \begin{equation*}
  \ee^{t D_F} G(u) := G(\nE(t,u)),
  \end{equation*}
  we have (set $ G=\Id $ and $ G=F $, respectively)
  \begin{equation*}
  F(\nE(t,u)) = F(\ee^{t D_F} u) = \ee^{t D_F} F(u).
  \end{equation*}
  In this formalism,~\eqref{quaxi} assumes a more `symmetric flavour',
  as in the linear case (see~\cite{auzingeretal18a}),
  \begin{equation*}
  \pdone\nE(t,u)
  = \th \big( F(\ee^{t D_F}u) + \ee^{t D_F} F(u) \big).
  \end{equation*}
  \rev{However, in the present context this formalism is
  of little practical use,
  and we stick to explicit, classical notation.}
}
\begin{equation} \label{quaxi}
\pdone\nE(t,u)
= \th \big( F(\nE(t,u)) + \pdtwo\nE(t,u) \cdot F(u) \big),
\end{equation}
which is valid due to~\eqref{fi},
motivates the definition of the {\em symmetrized defect}
\begin{equation} \label{sdnonlin}
\nDs(t,u) =
\pdone\nS(t,u) - \th \big( F(\nS(t,u)) + \pdtwo\nS(t,u) \cdot F(u) \big),
\end{equation}
satisfying $ \nDs(t,u) = \nDc(t,u) + \Order(t^{p+1}) $
(see Remark~\ref{rem:DcD}).

\begin{theorem}\label{th:GALS}
In terms of the symmetrized defect~\eqref{sdnonlin},
the local error has the integral representation
\begin{equation} \label{GALS-identity}
\nL(t,u_0) =
\int_0^t \pdtwo\nE(\tfrac{t-s}{2},\nS(s,\nE(\tfrac{t-s}{2},u_0)))
         \cdot \nDs(s,\nE(\tfrac{t-s}{2},u_0))\,\dd s.
\end{equation}
\end{theorem}
\begin{figure}[!ht]
\begin{center}
\quad\includegraphics[width=0.9\textwidth]{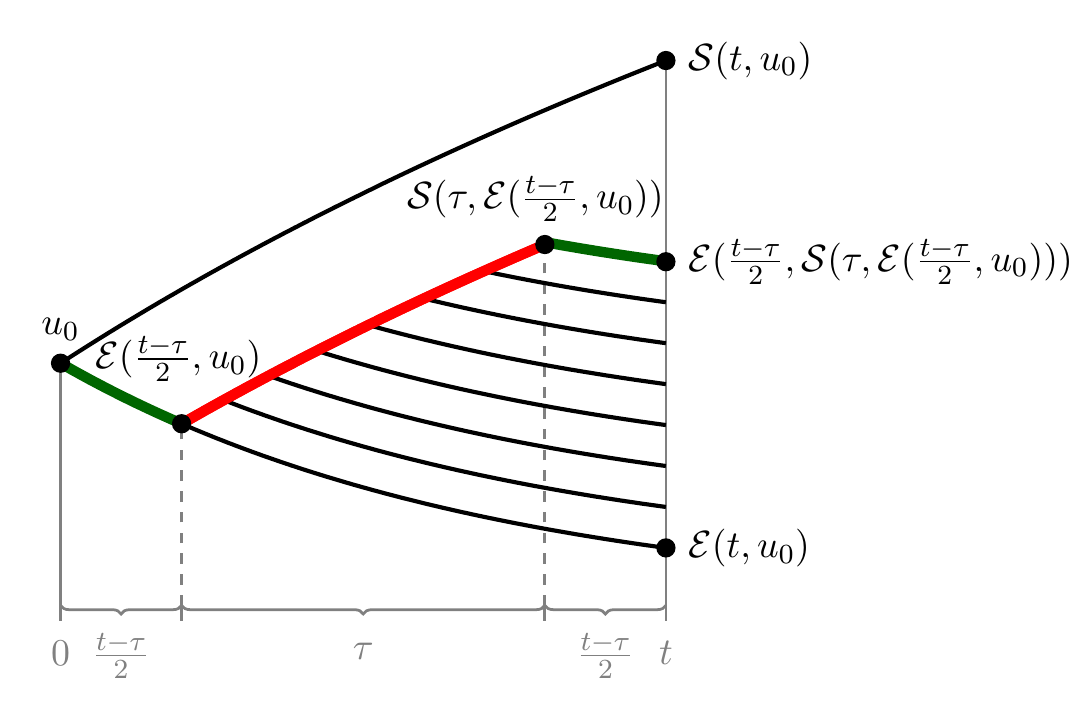}
\caption{Lady Windermere's Fan, Act 2\,\nicefrac{1}{2}
\label{fig:GALS}}
\end{center}
\end{figure}
{\em Proof.}
We reason in a similar way as in the proof of Theorem~\ref{th:GAL},
but now in the spirit of Figure~\ref{fig:GALS}.
For fixed $ t $, let
\begin{align*}
x(s) &= \nE(\tfrac{t-s}{2},u_0), \\
y(s) &= \nS(s,x(s)), \\
z(s) &= \nE(\tfrac{t-s}{2},y(s)).
\end{align*}
In this notation, we have
\begin{align*}
z(s) &= \nE(\tfrac{t-s}{2},\nS(s,\nE(\tfrac{t-s}{2},u_0))), \\
\text{satisfying} \quad z(0) &= \nE(t,u_0), ~~z(t) = \nS(t,u_0).
\end{align*}
Thus,
\begin{equation} \label{Lint-2}
\nL(t,u_0) = \nS(t,u_0) - \nE(t,u_0) = \int_0^t \ds z(s) \,\dd s,
\end{equation}
with
\begin{equation*}
\ds z(s)  = -\th F(z(s)) + \pdtwo\nE(\tfrac{t-s}{2},y(s)) \cdot \ds y(s).
\end{equation*}
Now, using~\eqref{fi}\footnote{{\em Mutatis mutandis:}
                               $ s,\frac{t-s}{2} $ and $ y(s) $
                               play the role of $ 0,t $ and $ u_0 $
                               from~\eqref{fi}.}
this can be rewritten in the form
\begin{subequations}
\begin{equation} \label{GALS-identity-proof-2}
\begin{aligned}
\ds z(s)
&= \th \ub{\big( -F(\nE(\tfrac{t-s}{2},y(s))) + \pdtwo\nE(\tfrac{t-s}{2},y(s)) \cdot F(y(s)) \big)}{=\,0} \\
& \quad {} + \pdtwo\nE(\tfrac{t-s}{2},y(s)) \cdot \big( \ul{\ds y(s) - \th F(y(s))} \big).
\end{aligned}
\end{equation}
Furthermore, from the definition~\eqref{sdnonlin}
of $ \nDs(s,u) $, with $ u=x(s) $ we obtain
\begin{equation} \label{GALS-identity-proof-3}
\begin{aligned}
&\ul{\ds y(s) - \th F(y(s))} \\
&= \pdone\nS(s,x(s)) +
   \pdtwo\nS(s,x(s)) \cdot \big( -\th F(x(s)) \big) - \th F(y(s)) \\
&= \pdone\nS(s,x(s)) -
   \th \big( F(y(s)) + \pdtwo\nS(s,x(s)) \cdot F(x(s)) \big) \\
&= \nDs(s,x(s)).
\end{aligned}
\end{equation}
\end{subequations}
After inserting~\eqref{GALS-identity-proof-3} into~\eqref{GALS-identity-proof-2},
together with~\eqref{Lint-2} we obtain~\eqref{GALS-identity}.
\hfill $ \square $

\section{Classical and symmetrized defect-based local error estimation}
\label{sec:locerrest}

\paragraph{Defect-based local error estimate}
The idea is due to~\cite{auzingeretal18a,auzingeretal13b}.
Let $ \nD(t,u)=\nDc(t,u) $ or $ \nDs(t,u) $, respectively,
and denote the integrands in~\eqref{GAL-identity} respectively~\eqref{GALS-identity},
generically by $ \Theta(s) $. Due to order~$ p $ we have
$ \nD(s,u) = \Order(s^p ) $ and $ \Theta(s) = \Order(s^p) $, whence
\begin{equation} \label{LDest}
\begin{aligned}
\nL(t,u_0) =
\int_0^t \Theta(s)\,\dd s
&\approx
\int_0^t \tfrac{s^p}{p!}\,\Theta^{(p)}(0)\,\dd s =
\tfrac{t^{p+1}}{(p+1)!}\,\Theta^{(p)}(0) \\
&\approx \tfrac{t}{p+1}\,\Theta(t) =
\tfrac{t}{p+1}\nD(t,u_0). 
\end{aligned}
\end{equation}
Here, `$ \approx $' means asymptotic approximation at the level
$ \Order(t^{p+2}) $.
\rev{
This approximation can be interpreted as an Hermite-type quadrature
of order $ p+1 $ for the local error integral,
where the quadrature error depends on
$ \tfrac{\partial^{p+1}}{\partial s^{p+1}}\nD(s,u_0) = \Order(1) $
due to $ \nD(s,u_0) = \Order(s^p) $, whence
\begin{equation*}
\nL(t,u_0) = \tfrac{t}{p+1}\,\nD(t,u_0) + \Order(t^{p+2})
\quad \text{for $ \nD = \nDc $\, or \,$ \nD = \nDs $.}
\end{equation*}
}
For a precise analysis of the resulting quadrature error based on its Peano representation
for the classical case in concrete applications,
see for instance~\cite{auzingeretal18a,auzingeretal13a,auzingeretal13b}.

Next we show that for the self-adjoint case and using the
symmetrized defect~\eqref{sdnonlin} we even
have\footnote{For the linear constant coefficient case
              see~\cite[Theorem~1]{auzingeretal18a}.}
\begin{equation*}
\nL(t,u_0) = \tfrac{t}{p+1} \nDs(t,u_0) + \Order(t^{p+3}).
\end{equation*}
To this end we consider the corrected scheme
\begin{equation} \label{Sdach-nonlin}
\nSD(t,u) = \nS(t,u) - \tfrac{t}{p+1} \nDs(t,u),
\end{equation}
and we show that it is of (global) order $ p+2 $.

\begin{theorem} \label{th:Sdach}
Consider a self-adjoint one-step scheme of (even) order $ p \geq 2 $,
represented by its flow $ \nS(t,u) $ satisfying~\eqref{Ssymm},
applied to an evolution equation~\eqref{y'=F(y)}.
Then the corrected scheme~\eqref{Sdach-nonlin} is almost self-adjoint,
i.e.,
\begin{subequations} \label{Sdach-results}
\begin{equation} \label{Sdach-almost-selfadjoint}
\nSD(-t,\nSD(t,u_0)) = u_0 + \Order(t^{2p+2}).
\end{equation}
Moreover, the local error $ \nLD(t,u) = \nSD(t,u) - \nE(t,u) $
of the corrected scheme satisfies
\begin{equation} \label{Sdach-highord-general-nonlin}
\nLD(t,u_0) = \Order(t^{p+3}),
\end{equation}
\end{subequations}
i.e., $ \nSD $  has even order $ p+2 $.
\end{theorem}
{\em Proof.}
We consider
\begin{align*}
\nSD(-t,\nSD(t,u_0))
&= \nS\big(\!-\!t,\nSD(t,u_0)) + \tfrac{t}{p+1}\,\nDs(-t,\nSD(t,u_0)) \\
&= \nS\big(\!-\!t,\nS(t,u_0) - \tfrac{t}{p+1}\,\nDs(t,u_0)\big) \\
& \quad {} + \tfrac{t}{p+1}\,\nDs\big(\!-\!t,\nS(t,u_0) - \tfrac{t}{p+1}\,\nDs(t,u_0) \big),
\end{align*}
apply Taylor expansion, and make use of the assumption that
$ \nS $ is self-adjoint, and the fact that $ t \nDs(t,u_0) = \Order(t^{p+1}) $:
\begin{align} \label{SdachSdach}
&\nSD(-t,\nSD(t,u_0))
 = \ub{\nS(-t,\nS(t,u_0))}{=\,u_0} \\
& \quad {} + \pdtwo\nS(-t,\nS(t,u_0))
             \,\cdot\,\big(\!-\!\tfrac{t}{p+1}\,\nDs(t,u_0)\big) + \Order(t^{2p+2}) \notag \\
& \quad {} + \tfrac{t}{p+1}\,\nDs(-t,\nS(t,u_0)) + \Order(t^{2p+2}) \notag \\
&= u_0 - \tfrac{t}{p+1}\,
   \ub{\big( \pdtwo\nS(-t,\nS(t,u_0)) \cdot \nDs(t,u_0) -
             \nDs(-t,\nS(t,u_0)) \big)}
      {\text{\small\bf\em critical term}} \notag
+\, \Order(t^{2p+2}).
\end{align}
Now we collect the contributions to the $ {\text{\bf\em critical term}} $.
First, from~\eqref{Ssymm}
we have\footnote{Here,
 $ \pdt\nS(-t,\nS(t,u_0)) $ means $ \pdt{\tilde\nS}(t,u_0) $
 with $ {\tilde\nS}(t,u_0) = \nS(-t,\nS(t,u_0)) $.}
\begin{equation*}
0 = \pdt\ub{\nS(-t,\nS(t,u_0))}{=\,u_0}
= -\pdone\nS(-t,\nS(t,u_0))
    + \pdtwo\nS(-t,\nS(t,u_0))\,\cdot\,\pdone\nS(t,u_0).
\end{equation*}
This implies
\begin{align*}
&\nDs(-t,\nS(t,u_0)) = \\
&= \pdone\nS(-t,\nS(t,u_0))
   - \th \big( F(\ub{\nS(-t,\nS(t,u_0))}{=\,u_0}) + \pdtwo\nS(-t,\nS(t,u_0))
      \cdot F(\nS(t,u_0)) \big) \\
&= \pdtwo\nS(-t,\nS(t,u_0)) \cdot \pdone\nS(t,u_0)
   - \th F(u_0) - \th \pdtwo\nS(-t,\nS(t,u_0)) \cdot F(\nS(t,u_0)) \\
&= \pdtwo\nS(-t,\nS(t,u_0)) \cdot
   \big( \pdone\nS(t,u_0) - \th F(\nS(t,u_0)) \big) - \th F(u_0).
\end{align*}
Summarizing and collecting terms gives
\begin{align*}
&\text{\normalsize\bf\em critical term} \,= \\
&= \pdtwo\nS(-t,\nS(t,u_0)) \cdot \nDs(t,u_0)
   - \nDs(-t,\nS(t,u_0)) \\
&= \pdtwo\nS(-t,\nS(t,u_0)) \cdot
     \Big( \canc{\,\pdone\nS(t,u_0)
             - \th  F(\nS(t,u_0))}
              - \th \pdtwo\nS(t,u_0) \cdot F(u_0) \Big) \\
& \quad {} -
  \pdtwo\nS(-t,\nS(t,u_0)) \cdot
   \Big( \canc{\,\pdone\nS(t,u_0) - \th F(\nS(t,u_0))}\,\Big)
          - \th F(u_0) \\
&= -\th \big( \pdtwo\nS(-t,\nS(t,u_0)) \cdot \pdtwo\nS(t,u_0)
               - \Id \big) \cdot F(u_0) \\
&= -\th \big( \ub{\pdunot\,\nS(-t,\nS(t,u_0))}{=\;\Id} - \,\Id \big) \cdot F(u_0) = 0.
\end{align*}
Thus,~\eqref{SdachSdach} indeed simplifies to~\eqref{Sdach-almost-selfadjoint},
\begin{equation*}
\nSD(-t,\nSD(t,u_0)) = u_0 + \Order(t^{2p+2}).
\end{equation*}
The proof of~\eqref{Sdach-highord-general-nonlin}
now works in the same way as for the linear case~\cite[proof of Theorem 1]{auzingeretal18a},
following the argument from~\cite[Theorem~II.3.2]{haireretal02}.
\hfill $ \square $

Assertion~\eqref{Sdach-highord-general-nonlin} is equivalent
to the fact that the symmetrized defect-based local error estimator
according to~\eqref{LDest},
\begin{equation} \label{improved-errest}
\nLLtildes(t,u_0) := \tfrac{t}{p+1} \nDs(t,u_0)
\end{equation}
is indeed of a better asymptotic quality than the classical defect,
with a deviation
\begin{equation} \label{p+3}
\nLLtildes(t,u_0) - \nL(t,u_0) = \Order(t^{p+3}),
\end{equation}
and not only $ \Order(t^{p+2}) $.

In the following sections we present some examples of self-adjoint
methods and show how to evaluate the symmetrized defect $ \nDs(t,u_0) $
as the basis for evaluating the local error estimator~\eqref{improved-errest}.

\section{Examples for the autonomous case} \label{sec:appl-auto}

\subsection{Example: Implicit midpoint rule} \label{subsec:imr}
We illustrate the defect computation for the simplest example
of a self-adjoint implicit one-step integrator.
The flow of the second order implicit midpoint rule
is defined by the relation
\begin{equation*}
\nS(t,u) = u + t\,F(\th(u+\nS(t,u))).
\end{equation*}
With
\begin{equation}\label{weh}
w = \nS(t,u)
\end{equation}
we obtain
\begin{subequations} \label{dS-imr}
\begin{equation*}
\pdone\nS(t,u)
= \ub{F(\th(u+w))}{=\,(w-u)/t}
  +\, t\,F'\big( \th(u+w) \big) \cdot \th \pdone\nS(t,u).
\end{equation*}
Thus, $ x = \pdone\nS(t,u) $ is obtained by solving the linear system
\begin{equation} \label{dS-imr-t}
\big( \Id - \tfrac{t}{2}F'(\th(u+w)) \big) \cdot x = F(\th(u+w)).
\end{equation}
Furthermore,
\begin{align*}
\pdtwo\nS(t,u)
&= \Id + t F'(\th(u+\nS(t,u))) \cdot \big(\th (\Id + \pdtwo\nS(t,u)) \big) \\
&= \Id + \tfrac{t}{2} F'(\th(u+w))
   + \tfrac{t}{2} F'(\th(u+w)) \cdot \pdtwo\nS(t,u),
\end{align*}
whence
\begin{equation*}
\big( \Id - \tfrac{t}{2}F'(\th(u+w)) \big) \cdot \pdtwo\nS(t,u)
= \big( \Id + \th F'(\th(u+w)) \big).
\end{equation*}
Thus, $ y = \pdtwo\nS(t,u) \cdot F(u) $ is obtained by solving the linear system
\begin{equation} \label{dS-imr-u}
\big( \Id - \tfrac{t}{2}F'(\th(u+w)) \big) \cdot y
= \big( \Id + \tfrac{t}{2} F'(\th(u+w)) \big) \cdot F(u),
\end{equation}
\end{subequations}
with the same matrix as in~\eqref{dS-imr-t}.

This gives the following defect representations.
\begin{itemize}
\item Classical defect:
\begin{equation*}
\nDc(t,u) = x - F(w),
\end{equation*}
where $ x = \pdone\nS(t,u) $ is the solution of~\eqref{dS-imr-t}
and with $w$ from~\eqref{weh}.
\item Symmetrized defect:
\begin{equation*}
\nDs(t,u) =
x - \th (F(w) + y),
\end{equation*}
where $ x = \pdone\nS(t,u) $ is the solution of~\eqref{dS-imr-t},
and $ y = \pdtwo\nS(t,u) \cdot F(u) $ is the solution
of~\eqref{dS-imr-u}. This can also be
written in the form
\begin{equation*}
\nDs(t,u) = z - \th F(w),
\end{equation*}
where $ z = x - \th y $ is the solution of
\begin{align*}
\big( \Id - \tfrac{t}{2}F'(\th(u+w)) \big) \cdot z
&= F(\th(u+w)) - \th F(u)
   - \tfrac{t}{4} F'(\th(u+w)) \cdot F(u).
\end{align*}
Thus, the computation of the symmetrized defect requires only one additional
evaluation of $F$ as compared to the classical version.
\end{itemize}

\subsection{Example: Strang splitting applied to a semilinear evolution equation}\label{subsec:strang}
We consider a semilinear problem of the form
\begin{equation*}
\dt u(t) = F(u(t)) = A u(t) + B(u(t)), \quad u(0) = u_0.
\end{equation*}
Denoting the flow of the nonlinear part by $ \nE_B(t,u) $,
the second order self-adjoint Strang splitting scheme is given by
\begin{equation*}
\nS(t,u) = \ee^{\tth A} \nE_B\big(t,\ee^{\tth A} u\big).
\end{equation*}
Let
\begin{equation*}
v_1 = \ee^{\tth A} u,
\quad
v_2 = \nE_B(t,v_1),
\quad
w = \ee^{\tth A} v_2 = \nS(t,u).
\end{equation*}
Then,
\begin{align*}
\pdone\nS(t,u)
&= \th A \nS(t,u)
    + \ee^{\tth A} \big( \pdone\nE_B(t,v_1)
                      + \pdtwo\nE_B(t,v_1)
                       (\th A v_1)
              \big) \\
&= \th A w
       + \ee^{\tth A} \big( B(v_2) + \th \pdtwo\nE_B(t,v_1) (A v_1)
              \big),
\end{align*}
and
\begin{equation*}
\pdtwo\nS(t,u) (\xi)
 = \ee^{\tth A} \pdtwo\nE_B(t,v_1) \big( \ee^{\tth A} \xi \big).
\end{equation*}
This gives the following defect representations.
\begin{itemize}
\item Classical defect:
\begin{align}
\nDc(t,u)
&= \pdone\nS(t,u) - F(\nS(t,u)) \notag \\
&=  \ee^{\tth A}
    \big( B(v_2)
            + \th \pdtwo\nE_B(t,v_1) \cdot (A v_1)
    \big)
     - \th A w - B(w). \label{dc-strang}
\end{align}
\item Symmetrized defect:
\begin{align} \label{ds-strang}
&\nDs(t,u)
= \pdone\nS(t,u) - \th \big( F(\nS(t,u))
    + \pdtwo\nS(t,u) \cdot F(u) \big) \notag \\
&= \th \canc{A w}
       + \ee^{\tth A}
              \big( B(v_2)
                      + \th \pdtwo\nE_B(t,v_1) \cdot (A v_1)
              \big) \notag \\
& \quad {} - \th \Big(
                   \canc{A w} + B(w) +
                   \ee^{\tth A}
                    \pdtwo\nE_B(t,v_1)
                     \cdot \big( \ee^{\tth A} (A u + B(u)) \big)
                 \Big) \notag \\
&= \ee^{\tth A} B(v_2)
    + \canc{\th \ee^{\tth A} \pdtwo\nE_B(t,v_1)
       (A v_1)} \notag \\
& \quad {}
    - \th B(w)
    - \canc{\th \ee^{\tth A} \pdtwo\nE_B(t,v_1) (A v_1)}
     - \th \ee^{\tth A}\pdtwo\nE_B(t,v_1)
        \big( \ee^{\tth A} B(u) \big) \notag \\[\jot]
&= \ee^{\tth A} \left( B(v_2)
    - \th \pdtwo\nE_B(t,v_1)
       \big( \ee^{\tth A} B(u) \big) \right)
      - \th B(w).
\end{align}
\end{itemize}
Thus,~\eqref{dc-strang} resp.~\eqref{ds-strang} require
one evaluation of $ \pdtwo\nE_B(t,v_1) \cdot (\,\cdot\,)$, and
either one or two evaluations of $ \ee^{\tth A} (\,\cdot\,) $,
respectively.

\subsection{Algorithmic realization for higher order splitting methods}\label{subsec:algsplit}

\begin{figure}[h!]
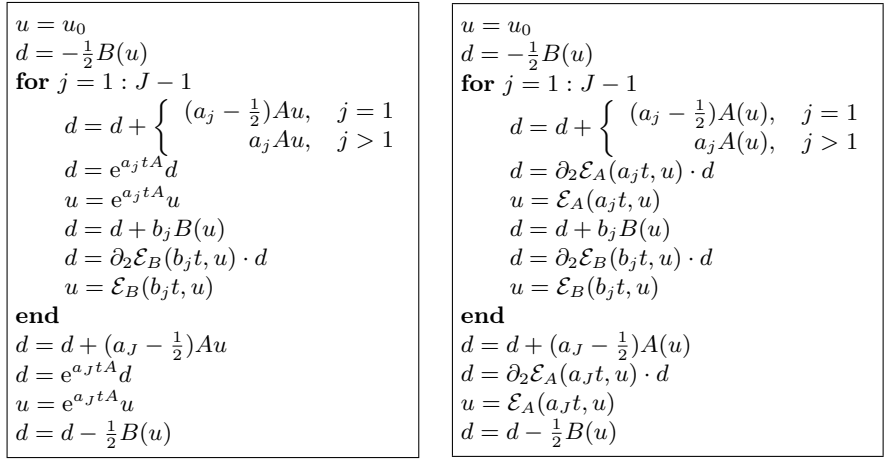

{\small
\begin{tabular}{cc}
&
\begin{boxedminipage}{0.72\textwidth}
\begin{tabbing}
 \qquad\= \kill
 $u=u_0$ \\
 $d=-\tfrac{1}{2}B(u)$ \\
 {\bf for} $j=1:J-1$ \\
 \>$d=d+\left\{
   \begin{array}{rl}
   (a_j-\tfrac{1}{2})Au,&j=1\\
   a_j Au,&j>1
   \end{array}
   \right.$\\
 \>$d=\ee^{a_j t A}d$\\
 \>$u=\ee^{a_j t A}u$\\
 \>$d=d+b_j B(u)$\\
 \>$d=\partial_2\nE_B(b_j t,u)\cdot d$\\
 \>$u=\nE_B(b_j t,u)$\\
 {\bf end}\\
 $d=d+(a_J-\tfrac{1}{2})Au$\\
 $d=\ee^{a_J t A}d$\\
 $u=\ee^{a_J t A}u$\\
 $d=d-\frac{1}{2}B(u)$
\end{tabbing}
\end{boxedminipage}
\quad
\begin{boxedminipage}{0.72\textwidth}
\begin{tabbing}
 \qquad\= \kill
 $u=u_0$ \\
 $d=-\tfrac{1}{2}B(u)$ \\
 {\bf for} $j=1:J-1$ \\
 \>$d=d+\left\{
   \begin{array}{rl}
   (a_j-\tfrac{1}{2})A(u),&j=1\\
   a_j A(u),&j>1
   \end{array}
   \right.$\\
 \>$d=\partial_2\nE_A(a_j t,u)\cdot d$\\
 \>$u=\nE_A(a_j t,u)$\\
 \>$d=d+b_j B(u)$\\
 \>$d=\partial_2\nE_B(b_j t,u)\cdot d$\\
 \>$u=\nE_B(b_j t,u)$\\
 {\bf end}\\
 $d=d+(a_J-\tfrac{1}{2})A(u)$\\
 $d=\partial_2\nE_A(a_J t,u)\cdot d$\\
 $u=\nE_A(a_J t,u)$\\
 $d=d-\frac{1}{2}B(u)$
\end{tabbing}
\end{boxedminipage}
\end{tabular}
}
\caption{Algorithmic realization of the symmetrized defect for splitting methods.\newline
Left: semilinear case. Right: nonlinear case.\label{fig:algssplit}}
\end{figure}

In Figure~\ref{fig:algssplit}, we give pseudocodes for the economical algorithmic
realization of the symmetrized defect when it is employed in the context
of splitting methods involving an arbitrary number of $J$ compositions.
If we denote the subflow of the nonlinear operator
by $\nE_B(t,u_0)$, an $n$-stage splitting approximation is defined by a composition of the two subflows,
\begin{equation*}
\nS(t,u_0) = \nE_B(b_J t, \cdots \ee^{a_2 t A}\nE_B(b_1 t, \ee^{a_1 t A}u_0)\cdots ).
\end{equation*}
An optimized fourth order method we will use in Section~\ref{subsec:nls} has the coefficient tableau given in Table~\ref{splitcoeffs}.
\begin{table}[ht]
\begin{center}
\begin{tabular}{|r||r|r|}
\hline
\multicolumn{1}{|c||}{$i$} & \multicolumn{1}{c|}{$a_i$} & \multicolumn{1}{c|}{$b_i$} \\ \hline\hline
1 &  0.267171359000977615 & $-$0.361837907604416033 \\
2 & $-$0.033827909669505667 &  0.861837907604416033 \\
3 &  0.533313101337056104 &  0.861837907604416033 \\
4 & $-$0.033827909669505667 & $-$0.361837907604416033 \\
5 &  0.267171359000977615 &  0 \\
\hline
\end{tabular}
\caption{Coefficients of the self-adjoint splitting method from~\cite[\texttt{Emb 4/3 AK s}]{splithp}.\label{splitcoeffs}}
\end{center}
\end{table}

The algorithms in Figure~\ref{fig:algssplit} have
the splitting approximation~$ u=\nS(t,u_0) $ and the symmetrized
defect $d=\nDs(t,u_0)$ as the output; for efficiency, $ u $ and $ d $ are
evaluated simultaneously. The left algorithm refers to the situation
where the operator $ A $ is linear, and
on the right the general nonlinear case is elaborated.

\section{The nonautonomous case, with examples} \label{sec:nonaut}
The results from Sections~\ref{sec:def} and~\ref{sec:locerrest}
carry over to nonautonomous evolution equations
\begin{subequations} \label{y'=F(t,y)}
\begin{equation} \label{y'=F(y)-t}
\dt u(t) = F(t,u(t)), \quad u(t_0)=u_0.
\end{equation}
For our purpose it is notationally more favorable to
introduce the~`local' variable~$ \tau $, such that
$ t = t_0 + \tau $, and reformulate~\eqref{y'=F(y)-t} in the form
\begin{equation} \label{y'=F(y)-ttau}
\dtau u(t_0+\tau) = F(t_0+\tau,u(t_0+\tau)), \quad u(t_0)=u_0.
\end{equation}
\end{subequations}
The exact flow associated with~\eqref{y'=F(t,y)}
is denoted by $ \nE(\tau,t_0,u) $.
It satisfies\footnote{Again, $ \pdone\nE(\tau,t_0,u_0) $ denotes
                      $ \dtau\nE(\tau,t_0,u_0) $,
                      and $ \pdtwo,\pdthree $ are defined analogously.}
\begin{equation*}
\pdone\nE(\tau,t_0,u_0) =
F(t_0+\tau,\nE(\tau,t_0,u_0)),
\quad \nE(0,t_0,u_0) = u_0.
\end{equation*}
\rev{
To infer the appropriate definition of the symmetrized
defect in this case there are two approaches, which we
both discuss for the sake of completeness. The first
one relies on a direct extension of the
fundamental identity~\eqref{fi} (Lemma~\ref{lem:fi}),
see Lemma~\ref{lem:fin} below. The other approach is based
on reformulating~\eqref{y'=F(t,y)} in autonomous form in the usual way,
leading to the same conclusion and showing that the theoretical
background based on Theorems~\ref{th:GALS} and~\ref{th:Sdach}
directly carries over to the nonautonomous case.
\begin{lemma} \label{lem:fin}
\begin{equation} \label{fin}
\begin{aligned}
&[\,\pdone\nE(\tau,t_0,u_0) =\,] \\
& \quad F(t_0+\tau,\nE(\tau,t_0,u_0))
= \pdtwo\nE(\tau,t_0,u_0) + \pdthree\nE(\tau,t_0,u_0) \cdot F(t_0,u_0).
\end{aligned}
\end{equation}
\end{lemma}
{\em Proof.}
The idea is the same as in the proof of Lemma~\ref{lem:fi}.
We proceed from the identity
\begin{equation*}
\nE(\tau+\sig,t_0,u_0) = \nE(\tau,t_0+\sig,\nE(\sig,t_0,u_0)).
\end{equation*}
Differentiation with respect to~$ \sig $ gives
\begin{align*}
\pdsig\nE(\tau+\sig,t_0,u_0)
&= \pdone\nE(\tau+\sig,t_0,u_0), \\
\pdsig\nE(\tau+\sig,t_0,u_0)\big|_{\sig=0}
&= \pdone\nE(\tau,t_0,u_0) = F(t_0+\tau,\nE(\tau,t_0,u_0)),
\end{align*}
and on the other hand,
\begin{align*}
&\pdsig\nE(\tau,t_0+\sig,\nE(\sig,t_0,u_0)) \\
&~~= \pdtwo\nE(\tau,t_0+\sig,\nE(\sig,t_0,u_0))
   + \pdthree\nE(\tau,t_0+\sig,\nE(\sig,t_0,u_0))
     \cdot \pdone\nE(\sig,t_0,u_0), \\
&\pdsig\nE(\tau,t_0+\sig,\nE(\sig,t_0,u_0))\big|_{\sig=0} \\
&~~= \pdtwo\nE(\tau,t_0,\nE(0,t_0,u_0))
   + \pdthree\nE(\tau,t_0,\nE(0,t_0,u_0))
     \cdot \pdone\nE(0,t_0,u_0) \\
&~~= \pdtwo\nE(\tau,t_0,u_0)
   + \pdthree\nE(\tau,t_0,u_0)
     \cdot F(t_0,u_0),
\end{align*}
which completes the proof.
\hfill $ \square $
}

\rev{
Alternatively, we can reformulate~\eqref{y'=F(y)-ttau} in autonomous form,
defining
\begin{equation*}
U = \left\lgroup \begin{array}{c}
      t_0+\tau \\
      u
    \end{array} \right\rgroup,
\quad
{\bm F}(U) = \left\lgroup \begin{array}{c}
       1 \\ F(t_0+\tau,u)
       \end{array} \right\rgroup
\end{equation*}
whence
\begin{equation*}
\dtau U(\tau) = {\bm F}(U(\tau)),
\quad
U(0) = \left\lgroup \begin{array}{c}
       t_0 \\ u_0
       \end{array} \right\rgroup,
\end{equation*}
and with the flow
\begin{equation*}
{\bm \nE}(\tau,U) =
{\bm \nE}(\tau,t_0,u) =
\left\lgroup \begin{array}{c}
 t_0+\tau \\ \nE(\tau,t_0,u)
\end{array} \right\rgroup
\end{equation*}
satisfying the fundamental identity according to Lemma~\ref{lem:fi},
\begin{equation} \label{fit}
[\,{\pdone{\bm \nE}}(\tau,U) =\,]~\;
{\bm F}({\bm \nE}(\tau,U)) = {\pdtwo{\bm \nE}}(\tau,U) \cdot {\bm F}(U).
\end{equation}
With $ U_0 = (t_0,u_0) $ we have
\begin{equation*} 
{\pdone{\bm \nE}}(\tau,U_0) = {\bm F}({\bm \nE}(\tau,U_0))
                     = \left\lgroup \begin{array}{c}
                        1 \\
                        F(t_0+\tau,\nE(\tau,t_0,u_0))
                       \end{array} \right\rgroup,
\quad {\bm \nE}(0,U_0) = U_0,
\end{equation*}
and
\begin{equation*}
{\pdtwo{\bm \nE}}(\tau,U_0) =
\left\lgroup \begin{array}{cc}
1 & 0 \\
\pdtwo\nE(\tau,t_0,u_0) &
\pdthree\nE(\tau,t_0,u_0)
\end{array} \right\rgroup.
\end{equation*}
Using~\eqref{fit} and evaluating the second component
again gives~\eqref{fin}.
}

\rev{
For a one-step approximation represented by
$ \nS(\tau,t_0,u_0) \approx \nE(\tau,t_0,u_0) $,
relation~\eqref{fin} again motivates the definition of the symmetrized defect
\begin{align}
&\nDs(\tau,t_0,u_0)
= \pdone\nS(\tau,t_0,u_0) \notag \\
& \quad {} - \th\big( F(t_0+\tau,\nS(\tau,t_0,u_0))
             + \pdtwo\nS(\tau,t_0,u_0)
              + \pdthree\nS(\tau,t_0,u_0)F(t_0,u_0)
                \big) \label{dst-n} \\
&= \big( \pdone - \th\pdtwo \big) \nS(\tau,t_0,u_0)
   - \th\big( F(t_0+\tau,\nS(\tau,t_0,u_0)) + \pdthree\nS(\tau,t_0,u_0)F(t_0,u_0) \big).
\notag
\end{align}
}

\paragraph{The linear nonautonomous case}
Now we consider the case of a linear time-dependent problem
\begin{equation} \label{y'=A(t)y}
\dtau u(t_0+\tau) = A(t_0+\tau) u(t_0+\tau), \quad u(t_0) = u_0.
\end{equation}
Since in the present case the flow
is linear in $ u_0 $, we write it in the simplified
form\footnote{\eqref{flowdet-magnus}
 is a minor abuse of notation.
 Note that $ \nE(\tau,t_0) $ can be expressed
 as a matrix exponential
 via the so-called Magnus expansion,
 see for instance~\cite{auzingeretal18a,blanesetal08b}.}
\begin{equation} \label{flowdet-magnus}
\nE(\tau,t_0,u_0) =: \nE(\tau,t_0) u_0,
\end{equation}
satisfying
\begin{equation*}
\pdone\nE(\tau,t_0) =
A(t_0+\tau)\,\nE(\tau,t_0),
\quad \nE(0,t_0) = \Id.
\end{equation*}
Note that
\begin{equation} \label{nonaut-timereversed}
\nE(-\tau,t_0+\tau)\,\nE(\tau,t_0) = \Id.
\end{equation}
A one-step approximation $ \nS(\tau,t_0,u_0) \approx \nE(\tau,t_0,u_0) $,
is also typically linear in~$ u_0 $,
\begin{equation*}
\nS(\tau,t_0,u_0) =: \nS(\tau,t_0) u_0 \approx \nE(\tau,t_0) u_0.
\end{equation*}
In particular, we again focus on self-adjoint schemes which are
characterized by the identity (cf.~\eqref{nonaut-timereversed})
\begin{equation} \label{nonaut-symmscheme}
\nS(-\tau,t_0+\tau)\,\nS(\tau,t_0) = \Id.
\end{equation}
For $ {\nS(\tau,t_0)u_0} $
we obtain the following defect representations.
\begin{itemize}
\item
Classical defect:
\begin{equation*}
\nDc(\tau,t_0,u_0) =: \nDc(\tau,t_0)u_0,
\end{equation*}
with
\begin{equation} \label{dct}
\nDc(\tau,t_0) =
\pdone\nS(\tau,t_0) - A(t_0+\tau) \nS(\tau,t_0).
\end{equation}
\item
Symmetrized defect~\eqref{dst-n}:
\begin{equation*}
\nDs(\tau,t_0,u_0) =: \nDs(\tau,t_0)u_0,
\end{equation*}
with
\begin{align}
&\nDs(\tau,t_0) =
\pdone\nS(\tau,t_0) -
\th\big( A(t_0+\tau) \nS(\tau,t_0)
          + \pdtwo\nS(\tau,t_0)
           + \nS(\tau,t_0)A(t_0)
   \big) \notag \\
& \quad = \big( \pdone - \th\pdtwo \big) \nS(\tau,t_0)
- \th\big( A(t_0+\tau) \nS(\tau,t_0) + \nS(\tau,t_0) A(t_0) \big). \label{dst}
\end{align}
\end{itemize}

\subsection{Example: Exponential midpoint rule}\label{subsec:emr}
The self-adjoint second order exponential midpoint rule
applied to~\eqref{y'=A(t)y} is given by
\begin{equation*} 
S(\tau,t_0) = \ee^{\tau A(t_0+\frac{\tau}{2})}.
\end{equation*}
Let
\begin{equation*}
\nR(\tau,t_0)(\,\cdot\,) =
\tfrac{\dd}{\dd\Omega} \ee^{\Omega}\big|_{\Omega=\tau A(t_0+\frac{\tau}{2})}(\,\cdot\,),
\end{equation*}
where $ \tfrac{\dd}{\dd\Omega} \ee^\Omega $ denotes the Fr{\'e}chet
derivative of the matrix exponential, see~\eqref{mimi} below. Then,
\begin{align*}
\pdone \nS(\tau,t_0)
&= \nR(\tau,t_0) \big( \pdtau(\tau A(t_0+\tfrac{\tau}{2}) \big) \\
&= \nR(\tau,t_0) \big( A(t_0+\tfrac{\tau}{2})
                     + \th\tau A'(t_0+\tfrac{\tau}{2}) \big), \\
\pdtwo \nS(\tau,t_0)
&= \nR(\tau,t_0) \big( \pdtnot(\tau A(t_0+\tfrac{\tau}{2}) \big) \\
&= \nR(\tau,t_0) \big( \tau A'(t_0+\tfrac{\tau}{2}) \big).
\end{align*}
This gives the following defect representations.
\begin{itemize}
\item
Classical defect~\eqref{dct}:
\begin{equation} \label{dcA(t)}
\nDc(\tau,t_0) =
\nR(\tau,t_0) \big( A(t_0+\tfrac{\tau}{2})
                    + \th\tau A'(t_0+\tfrac{\tau}{2}) \big)
- A(t_0+\tau) \nS(\tau,t_0).
\end{equation}
\item
Symmetrized defect~\eqref{dst}:
\begin{align}
&\nDs(\tau,t_0) =
\nR(\tau,t_0) \big( A(t_0+\tfrac{\tau}{2})
                     + \canc{\th\tau A'(t_0+\tfrac{\tau}{2})} \big) \notag \\
& \qquad\qquad {} - \th\big( A(t_0+\tau) \nS(\tau,t_0)
                      + \nR(\tau,t_0) \big( \canc{\tau A'(t_0+\tfrac{\tau}{2})}
                       - \nS(\tau,t_0) A(t_0)
                \big) \notag \\
& \quad {} = \nR(\tau,t_0) \big( A(t_0+\tfrac{\tau}{2}) \big)
   - \th \big( A(t_0+\tau) \nS(\tau,t_0) + \nS(\tau,t_0) A(t_0) \big). \label{dsA(t)}
\end{align}
\end{itemize}
Here, the explicit representation
\begin{align} \label{mimi}
\nR(\tau,t_0) \big( V \big)
&= \int_0^1 \ee^{\sig\tau A(t_0+\frac{\tau}{2})}
            V
            \ee^{(1-\sig)\tau A(t_0+\frac{\tau}{2})}\,\dd\sig \notag \\
&= \int_0^1 \ee^{\sig\tau A(t_0+\frac{\tau}{2})}
            V
            \ee^{-\sig\tau A(t_0+\frac{\tau}{2})}\,\dd\sig
   \cdot \nS(\tau,t_0)
\end{align}
follows from~\cite[(10.15)]{higham08}.
For evaluating~\eqref{dcA(t)}, a sufficiently accurate quadrature
approximation for the integral according to~\eqref{mimi} is required.
This involves evaluation of $ A' $ and the commutator $ [A,A'] $,
see~\cite{auzingeretal18a}.
In contrast, the relevant term from~\eqref{dsA(t)} simplifies to
\begin{align*}
\nR(\tau,t_0) \big( A(t_0+\tfrac{\tau}{2}) \big)
&=
\int_0^1 \ee^{\sig\tau A(t_0+\frac{\tau}{2})}
         A(t_0+\tfrac{\tau}{2})\,
         \ee^{-\sig\tau A(t_0+\frac{\tau}{2})}\,\dd\sig
   \cdot \nS(\tau,t_0) \\
&= A(t_0+\tfrac{\tau}{2}) \nS(\tau,t_0)
 = \nS(\tau,t_0) A(t_0+\tfrac{\tau}{2}),
\end{align*}
whence the symmetrized defect~\eqref{dsA(t)} can be evaluated exactly,
\begin{equation} \label{dsA(t)-1}
\begin{aligned}
\nDs(\tau,t_0)
&= \big( A(t_0+\tfrac{\tau}{2}) - \th A(t_0+\tau) \big) \nS(\tau,t_0)
   - \th \nS(\tau,t_0) A(t_0) \\
&= \nS(\tau,t_0) \big( A(t_0+\tfrac{\tau}{2}) - \th A(t_0) \big)
  - \th A(t_0+\tau) \nS(\tau,t_0).
\end{aligned}
\end{equation}
This involves an additional application of $ \nS(\tau,t_0) $,
but it does not require evaluation of the derivative $ A' $
or of a commutator expression.
We also note that the applications of
$\nS$ from left and right can be evaluated in parallel.

\subsection{Algorithmic realization for higher order Magnus-type methods} \label{subsec:cfm}

The integrators which we consider for the
numerical approximation of~\eqref{zener}
are commutator-free Magnus-type methods (CFM) and classical Magnus integrators.

In contrast to the special case of the exponential midpoint rule,
for practical evaluation the defect needs to be approximated
in an asymptotically correct way.
To this end we require an approximation scheme
which preserves the desired order $ p+2 $ of the corrected
scheme~\eqref{Sdach-nonlin}, or equivalently,
the asymptotic quality~\eqref{p+3}
of the local error estimator is not affected by such an approximation.

Various versions of the resulting classical defect-based error estimators for
these exponential integrators are presented in~\cite{auzingeretal18b}.
We now follow two of these approaches.
To keep the presentation self-contained within reason,
we briefly recapitulate the underlying material
from~\cite[Section~3]{auzingeretal18b},
and we introduce the corresponding symmetrized defect approximations.

\subsubsection{Commutator-free Magnus-type integrators}
As the basic integrator we consider a
\emph{commutator-free Magnus-type (CFM)} method~\cite{alvfeh11},
\begin{subequations} \label{cfm}
\begin{equation} \label{cfmS}
\nS(\tau,t_0) = \nS_J(\tau,t_0)\,\cdots\,S_1(\tau,t_0),
\end{equation}
where
\begin{equation} \label{cfmSB}
\begin{aligned}
&S_j(\tau,t_0) = \ee^{\Omega_j(\tau,t_0)} = \ee^{\tau B_j(\tau,t_0)}, \\
&\text{with}~~ B_j(\tau,t_0) = \sum_{k=1}^{K} a_{jk}\,A(t_0+c_k\tau),
\end{aligned}
\end{equation}
\end{subequations}
where the coefficients $ c_k $ and $ a_{jk} $
are chosen in such a way that
a desired order of consistency is obtained.
Note that the assumption of symmetry of the scheme also implies symmetry of the
coefficients in the following sense,
\begin{subequations}\label{symmetriesgack}
\begin{equation}
c_k-\tfrac{1}{2} = \tfrac{1}{2}-c_{K+1-k},
\quad k=1,\ldots,K,
\end{equation}
and
\begin{equation}
a_{jk} = a_{J+1-j,K+1-k}, \quad j=1,\ldots,J, ~~ k=1,\ldots,K.
\end{equation}
\end{subequations}

Our construction involves evaluation of the derivatives
\begin{equation*}
\pdtau\ee^{\Omega_j(\tau,t_0)}
= \Gamma_{\tau,j}(\tau,t_0)\,\ee^{\Omega_j(\tau,t_0)},\quad
\tfrac{\partial}{\partial t_0}\ee^{\Omega_j(\tau,t_0)}
= \Gamma_{t_0,j}(\tau,t_0)\,\ee^{\Omega_j(\tau,t_0)},
\end{equation*}
where
\begin{equation*}
\Gamma_{\tau,j}(\tau,t_0)
= B_j(\tau,t_0) +
\sum_{m \geq 0}\tfrac{1}{(m+1)!}
 \tau^{m+1}\mathrm{ad}^m_{B_j(\tau,t_0)}
  (\pdtau B_j(\tau,t_0)),
\end{equation*}
and
\begin{equation*}
\Gamma_{t_0,j}(\tau,t_0)
= \sum_{m \geq 0}\tfrac{1}{(m+1)!}
   \tau^{m+1}\mathrm{ad}^m_{B_j(\tau,t_0)}
    (\pdtwo B_j(\tau,t_0)).
\end{equation*}
{Applying the product rule
to $ \nS(\tau,t_0) $ defined in~\eqref{cfm}
we see that the symmetrized} defect~\eqref{dst}
of the numerical approximation
is an expression involving the derivatives
\begin{equation} \label{wadlstrumpf}
\big( \pdone-\tfrac{1}{2}\,\pdtwo \big) \nS_j(\tau,t_0)
= \Gamma_{j}(\tau,t_0)\,\nS_j(\tau,t_0),
\end{equation}
with
\begin{align}
\Gamma_j(\tau,t_0)
&=\Gamma_{\tau,j}(\tau,t_0)-\tfrac{1}{2}\Gamma_{t_0,j}(\tau,t_0) \label{Gammaj-series} \\
&=B_j(\tau,t_0)
  +\sum_{m \geq 0}\tfrac{1}{(m+1)!}
   \tau^{m+1}\mathrm{ad}^m_{B_j(\tau,t_0)}
    ( \Bcheck_j(\tau,t_0)), \notag
\end{align}
where we have defined
\begin{equation*}
\Bcheck_j(\tau,t_0)
= \big( \pdone - \tfrac{1}{2} \pdtwo \big) B_j(\tau,t_0) =
\sum_{k=1}^{K} {a_{jk}(c_k-\tfrac{1}{2})}A'(t_0+c_k\tau).
\end{equation*}
One possible computable approximation is obtained by truncating the series (\ref{Gammaj-series});
we will refer to the resulting procedure as \emph{Taylor variant}.
The procedure in conjunction with the classical defect is given in detail
in~\cite[Section~3]{auzingeretal18b}.

We remark at this point that symmetry of the basic CFM integrator
implies that truncation of the series~\eqref{Gammaj-series} at $m=p$,
i.e., approximating $ \Gamma_j(\tau,t_0) $
by\footnote{A~priori one would expect that it is required to include
                  the term of degree~$ p+1 $ also.}
\begin{equation} \label{Gammaj-series-truncated}
{\widetilde\Gamma}_j(\tau,t_0)
= B_j(\tau,t_0)
  + \sum_{m=0}^{p-1} \tfrac{1}{(m+1)!}
     \tau^{m+1}\mathrm{ad}^m_{B_j(\tau,t_0)}
      (\Bcheck_j(\tau,t_0))
\end{equation}
is already sufficient to obtain a defect approximation of accuracy $p+2$,
as is demonstrated in the following.

\begin{figure}[ht!]
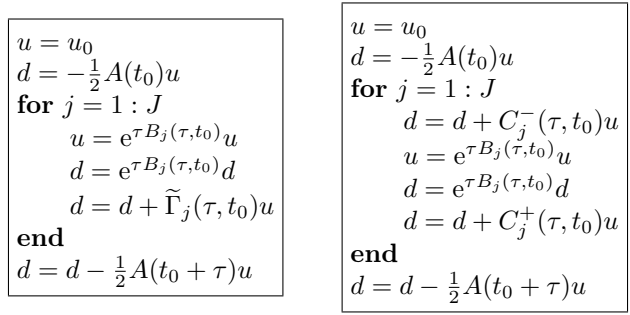

\begin{center}
\begin{tabular}{cc}
\begin{boxedminipage}{0.72\textwidth}
\begin{tabbing}
 \qquad\= \kill
 $u=u_0$ \\
 $d=-\tfrac{1}{2}A(t_0)u$ \\
 {\bf for} $j=1:J$\\
 \>$u=\ee^{\tau B_j(\tau,t_0)}u$\\
 \>$d=\ee^{\tau B_j(\tau,t_0)}d$\\
 \>$d=d+\widetilde{\Gamma}_j(\tau,t_0)u$\\
 {\bf end}\\
 $d = d-\tfrac{1}{2}A(t_0+\tau)u$
\end{tabbing}
\end{boxedminipage}\quad
&
\quad \begin{boxedminipage}{0.72\textwidth}
\begin{tabbing}
 \qquad\=\qquad\= \kill
 $u=u_0$ \\
 $d=-\tfrac{1}{2}A(t_0)u$\\
 {\bf for} $j=1:J$\\
 \>$d=d+C^{-}_j(\tau, t_0)u$\\
 \>$u=\ee^{\tau B_j(\tau,t_0)}u$\\
 \>$d=\ee^{\tau B_j(\tau,t_0)}d$\\
 \>$d=d+C^{+}_j(\tau, t_0)u$\\
 {\bf end}\\
 $d=d-\tfrac{1}{2}A(t_0+\tau)u$
\end{tabbing}
\end{boxedminipage}
\end{tabular}
\caption{Algorithmic realization of the symmetrized defect for CFM methods.\newline
Left:~Taylor variant.
Right:~Hermite variant. \label{fig:algscfm}}
\end{center}
\end{figure}

\begin{proposition}\label{prop:symm}
Let $\nDs$ be the symmetrized defect of a self-adjoint CFM integrator of
order $p$, and $\nDstilde$ its approximation constructed via the
truncated Taylor variant
according to~\eqref{Gammaj-series-truncated}. Then,
\begin{equation*}
\nDs(\tau,u_0) - \nDstilde(\tau,u_0) = \Order(\tau^{p+2}).
\end{equation*}
\end{proposition}
{\em Proof.}
Observe that
\begin{equation*}
B_j(\tau,t_0) = X_j A(t_0) + \Order(\tau),
\quad
\Bcheck_j(\tau,t_0)
= Y_j A'(t_0) + \Order(\tau),
\end{equation*}
where
\begin{equation*}
X_j = \sum_{k=1}^K a_{jk}, \quad
Y_j = \sum_{k=1}^K a_{jk}(c_k-\tfrac{1}{2}).
\end{equation*}
Thus,
\begin{equation*}
\Gamma_j(\tau,t_0)-{\widetilde\Gamma}_j(\tau,t_0) =
\tfrac{1}{(p+1)!}\tau^{p+1} X_j^p Y_j\,\mathrm{ad}_{A(t_0)}^{p}(A'(t_0))
+ \Order(\tau^{p+2}).
\end{equation*}
Inserting this in the computational algorithm given in Figure~\ref{fig:algscfm} (left)
and taking into account that
\begin{equation*}
\ee^{\tau B_j(\tau,t_0)} = \mathrm{Id} + \Order(\tau),
\end{equation*}
the total error resulting from substitution of the exact defect $\nDs$ by the
truncated Taylor approximation of $\Gamma_j$ is
\begin{equation*}
\nDs(\tau,t_0)-\nDstilde(\tau,t_0) = \tfrac{1}{(p+1)!} \tau^{p+1} Z\,
\mathrm{ad}_{A(t_0)}^{p}(A'(t_0)) + \Order(\tau^{p+2}),
\end{equation*}
with
\begin{equation*}
Z = \sum_{j=1}^{J} X_j^p \, Y_j.
\end{equation*}
To establish the assertion of the proposition we now show $Z=0$:
From~\eqref{symmetriesgack},
\begin{equation*}
X_j^p\,Y_j = -X_{J+1-j}^p Y_{J+1-j}, ~~j=1,\ldots,J,
\quad X_{\lfloor J/2 \rfloor +1}^p Y_{\lfloor J/2 \rfloor +1}
= 0 \ \text{ if } J \text{ is odd},
\end{equation*}
whence
\begin{equation*}
Z = \sum_j^{\lfloor J/2 \rfloor} (X_j^p Y_j + X_{J+1-j}^p Y_{J+1-j})
~~\big[+ X_{\lfloor J/2 \rfloor +1}^p Y_{\lfloor J/2 \rfloor +1}
\ \text{ if }
J \text{ is odd} \big] = 0,
\end{equation*}
which completes the proof. \hfill $ \square $

As an alternative to the series representation~\eqref{Gammaj-series},
we may use the integral {representation} which follows from~\cite[(10.15)]{higham08},
\begin{equation*} \label{Gammaj-int}
\Gamma_j(\tau,t_0) = B_j(\tau,t_0)+
\int_0^\tau \ee^{\sig B_{j}(\tau,t_0)} \Bcheck_j(\tau,t_0)\,\ee^{-\sig B_{j}(\tau,t_0)}\, \dd\sig,
\end{equation*}
and apply a {$ p $-th order two-sided Hermite-type quadrature}
(see~\cite[Section~3]{auzingeretal18b})
to approximate the integral.
We will refer to the resulting procedure as \emph{Hermite variant}.
The procedure in conjunction with the classical defect
was also introduced in~\cite[Section~3]{auzingeretal18b}.
Similarly as for the Taylor variant,
it can be shown that quadrature of order $p$ is sufficient
to obtain a defect approximation of order $p+2$.

These two sketched strategies result in the
procedures given as pseudocode in Figure~\ref{fig:algscfm} where
the defect $d=\nDs(\tau,t_0)u_0$ is computed as the output
{along with the basic approximation $ u = \nS(\tau,t_0)u_0 $}.
Then, for order $p=4$, for instance, for the Taylor variant we have
\begin{align*}
{\widetilde\Gamma}_j(\tau,t_0)
&= B_j(\tau,t_0) + \tau\Bcheck_j(\tau,t_0)
+\tfrac{1}{2}\tau^2[B_j(\tau,t_0),\Bcheck_j(\tau,t_0)] \notag \\
& \quad {} +\tfrac{1}{6}\tau^3[B_j(\tau,t_0),[B_j(\tau,t_0),\Bcheck_j(\tau,t_0)]] \\
& \quad {} +\tfrac{1}{24}\tau^4[B_j(\tau,t_0),[B_j(\tau,t_0),[B_j(\tau,t_0),\Bcheck_j(\tau,t_0)]]],
\end{align*}
and for the Hermite variant,
\begin{equation*}
C^{\pm}_j(\tau,t_0)=\tfrac{1}{2}\big(B_j(\tau,t_0)+\tau\Bcheck_j(\tau,t_0)\big)
\pm\tfrac{1}{12}\tau^2[B_j(\tau,t_0),\Bcheck_j(\tau,t_0)].
\end{equation*}

\subsubsection{Classical Magnus integrators}

As an example we consider the classical fourth order
Magnus integrator based on quadrature at Gaussian points
(see~\cite{auzingeretal18b}),
\begin{subequations} \label{MC-4}
\begin{equation} \label{MC-4a}
\nS(\tau,t_0) = \ee^{\Omega(\tau,t_0)} = \ee^{\tau B(\tau,t_0)},
\end{equation}
where $ {\Omega(\tau,t_0)} = \tau B(\tau,t_0)$
approximates the Magnus series $ {\bm\Omega}(\tau,t_0) $,
\begin{align}
B(\tau,t_0)
&= \tfrac{1}{2}\big(A(t_0+c_1\tau)+A(t_0+c_2\tau)\big)
-\tfrac{\sqrt{3}}{12}\tau[A(t_0+c_1\tau),A(t_0+c_2\tau)], \notag \\
c_{1,2}&=\tfrac{1}{2}\pm\tfrac{\sqrt{3}}{6}. \label{MC4-b}
\end{align}
\end{subequations}
Following~\cite[Section~3]{auzingeretal18b} for
the classical defect, the symmetrized defect~\eqref{dst} is now given by
\begin{equation*}
\nDs(\tau,t_0)
= \big(\Gamma(\tau,t_0) -
\tfrac{1}{2}A(t_0+\tau)\big)\,\nS(\tau,t_0)-\tfrac{1}{2}\nS(\tau,t_0)A(t_0),
\end{equation*}
where $ \Gamma(\tau,t_0) $ has a series representation
analogous to~\eqref{Gammaj-series}.
To approximate $ \nDs(\tau,t_0) $ in an asymptotically correct way,
we again truncate the series defining $ \Gamma(\tau,t_0) $ and obtain
the \emph{Taylor variant}
\begin{equation*}
\nDs(\tau,t_0)
\approx \big(\widetilde{\Gamma}(\tau,t_0)
  -\tfrac{1}{2}A(t_0+\tau)\big)\nS(\tau,t_0)-\tfrac{1}{2}\nS(\tau,t_0)A(t_0),
\end{equation*}
where
\begin{align*}
{\widetilde{\Gamma}}(\tau,t_0)
&= B(\tau,t_0) + \tau\Bcheck(\tau,t_0)
+\tfrac{1}{2}\tau^2[B(\tau,t_0),\Bcheck_j(\tau,t_0)] \\
& \quad {} +\tfrac{1}{6}\tau^3[B(\tau,t_0),[B(\tau,t_0),\Bcheck(\tau,t_0)]] \\
& \quad {} +\tfrac{1}{24}\tau^4[B(\tau,t_0),[B(\tau,t_0),[B(\tau,t_0),\Bcheck(\tau,t_0)]]],
\end{align*}
with
\begin{equation} \label{tildeBMagnus}
\begin{aligned}
\Bcheck(\tau,t_0)
&= \big( \pdone - \tfrac{1}{2} \pdtwo \big)B(\tau,t_0) \\
&= \tfrac{1}{2}\big((c_1-\tfrac{1}{2})A'(t_0+c_1\tau)+(c_2-\tfrac{1}{2})A'(t_0+c_2\tau)\big) \\
& \quad {} -\tfrac{\sqrt{3}}{12}[A(t_0+c_1\tau),A(t_0+c_2\tau)] \\
& \quad {} -\tfrac{\sqrt{3}}{12}(c_1-\tfrac{1}{2})\tau\,[A'(t_0+{c_1}\tau),A(t_0+c_2\tau)] \\
& \quad {} -\tfrac{\sqrt{3}}{12}(c_2-\tfrac{1}{2})\tau\,[A(t_0+c_1\tau),A'(t_0+c_2\tau)].
\end{aligned}
\end{equation}
Due to $c_1+c_2=1$ it follows by expansion in $\tau$ that
$\Bcheck(\tau,t_0)=\Order(\tau)$. Thus, truncation after $p=4$ again yields
a sufficiently accurate approximation.
Alternatively, application of fourth order two-sided Hermite quadrature
for the approximation of $\Gamma(\tau,t_0)$ yields the \emph{Hermite variant}
\begin{equation*}
\nDs(\tau,t_0) \approx
\big(C^{+}(\tau,t_0)-\tfrac{1}{2}A(t_0+\tau)\big)\nS(\tau,t_0)+\nS(\tau,t_0)\big(C^{-}(\tau,t_0)-\tfrac{1}{2}A(t_0)\big),
\end{equation*}
where
\begin{equation*}
C^{\pm}(\tau,t_0)=\tfrac{1}{2}\big(B(\tau,t_0)+\tau\Bcheck(\tau,t_0)\big)
\pm\tfrac{1}{12}\tau^2[B(\tau,t_0),\Bcheck(\tau,t_0)],
\end{equation*}
with $\Bcheck(\tau,t_0)$ as in~\eqref{tildeBMagnus}.

\section{Numerical examples} \label{sec:num}
We illustrate the theoretical analysis of the deviation of the symmetrized
error estimator by showing the orders \rev{of the error} of the basic
integrator and of the deviation of the error estimator
from \rev{the true error}. We will consider
splitting methods for a cubic nonlinear Schr{\"o}dinger equation
and commutator-free and classical Magnus-type integrators for a Rosen--Zener model.

\subsection{Cubic Schr{\"o}dinger equation}\label{subsec:nls}

We solve the cubic nonlinear Schr{\"o}dinger equation on the real line $ x \in \RR $
\begin{equation} \label{eq:cnls}
\begin{aligned}
\ii\,\partial_t \psi(x,t)
&= -\,\tfrac{1}{2}\partial_x^2\,\psi(x,t) -|\psi(x,t)|^2\,\psi(x,t), \quad t>0, \\
\psi(x,0) &= \psi_0(x)
\end{aligned}
\end{equation}
by splitting methods. Here, a soliton solution exists,
\begin{equation*} 
\psi(x,t) = 2\,\ee^{\ii(\frac{3}{2}t-x)}\,\text{sech}(2(t+x))
\end{equation*}
Our initial condition is chosen commensurate with this solution,
and we truncate the spatial domain to $ x \in [-16,16] $ and impose periodic boundary
conditions. Spectral collocation at 512 equidistant mesh points leads to an ODE system of the form
\begin{equation*}
\dt \Psi(t) = F(\Psi(t)) = A \Psi(t) + B(\Psi(t)), \quad \Psi(0) = \Psi_0,
\end{equation*}
with $ A \Psi \sim \tfrac{\ii}{2} \partial_x^2\,\psi $ and $ B(\Psi) \sim \ii\, |\psi|^2 \psi. $
We solve this by the second order Strang splitting
and by the self-adjoint fourth-order method
represented by the higher-order method in the embedded pair referred to as \texttt{Emb 4/3 AK s}
in the collection~\cite{splithp}, recapitulated for easy reference in Table~\ref{splitcoeffs}
in Section~\ref{subsec:algsplit}.
The $ A $-part is solved via [I]FFT, while the $ B $-part can
be integrated directly on the given mesh.

In Table~\ref{tab:conv1}, we give the local error of the Strang splitting and the
error of our symmetrized error
estimator as compared to the exact errors.
\rev{Table~\ref{tab:conv1b} shows the global errors on the interval $[0,1/8]$ of the basic integrator
and of the solution corrected by adding the error estimate.
In accordance with our theory, we observe local orders three and five, respectively,
and the expected orders two and four for the global errors. Likewise, Table~\ref{tab:conv2} shows
orders five and seven for the local errors of the fourth order integrator from~\cite[\texttt{Emb 4/3 AK s}]{splithp},
and Table~\ref{tab:conv2b} shows the matching global errors}.

\begin{table}[!ht]
\begin{center}
\begin{tabular}{||r||r|r||c|r||}
\hline
\multicolumn{1}{||c||}{$\tau$}    & \multicolumn{1}{c|}{$\|\nL(\tau,u_0)\|_2$} & \multicolumn{1}{|c||}{order}   &
\multicolumn{1}{c|}{$\|{\nLLtildes}(\tau,u_0) - \nL(\tau,u_0)\|_2$} &  \multicolumn{1}{|c||}{order}   \\ \hline
1.563e$-$02 & 3.791e$-$05 & 2.98 & 3.377e$-$07 & 4.59 \\
7.813e$-$03 & 4.753e$-$06 & 3.00 & 1.161e$-$08 & 4.86 \\
3.906e$-$03 & 5.946e$-$07 & 3.00 & 3.726e$-$10 & 4.96 \\
1.953e$-$03 & 7.434e$-$08 & 3.00 & 1.172e$-$11 & 4.99 \\
9.766e$-$04 & 9.293e$-$09 & 3.00 & 3.669e$-$13 & 5.00 \\
4.883e$-$04 & 1.162e$-$09 & 3.00 & 1.160e$-$14 & 4.98 \\
\hline
\end{tabular}
\caption{Local error and deviation of the symmetrized defect-based error
estimator for the second order Strang splitting applied to~\eqref{eq:cnls}.\label{tab:conv1}}
\end{center}
\end{table}

\begin{table}[!ht]
\begin{center}
\begin{tabular}{||r||r|r||c|r||}
\hline
\multicolumn{1}{||c||}{$\tau$}    & \multicolumn{1}{c|}{global error} & \multicolumn{1}{|c||}{order}   &
\multicolumn{1}{c|}{error of corrected solution} &  \multicolumn{1}{|c||}{order}   \\ \hline
1.563e$-$02 & 2.539e$-$04 & 1.99 & 5.703e$-$07 & 4.00 \\
7.813e$-$03 & 6.354e$-$05 & 2.00 & 3.634e$-$08 & 3.97 \\
3.906e$-$03 & 1.589e$-$05 & 2.00 & 2.283e$-$09 & 3.99 \\
1.953e$-$03 & 3.972e$-$06 & 2.00 & 1.428e$-$10 & 4.00 \\
9.766e$-$04 & 9.931e$-$07 & 2.00 & 8.928e$-$12 & 4.00 \\
4.883e$-$04 & 2.483e$-$07 & 2.00 & 5.611e$-$13 & 3.99 \\
\hline
\end{tabular}
\caption{\rev{Global error and corrected solution for the second order Strang splitting applied to~\eqref{eq:cnls}.}\label{tab:conv1b}}
\end{center}
\end{table}

\begin{table}[!ht]
\begin{center}
\begin{tabular}{||r||r|r||c|r||}
\hline
\multicolumn{1}{||c||}{$\tau$}    & \multicolumn{1}{c|}{$\|\nL(\tau,u_0)\|_2$} & \multicolumn{1}{|c||}{order}   &
\multicolumn{1}{c|}{$\|{\nLLtildes}(\tau,u_0) - \nL(\tau,u_0)\|_2$} &  \multicolumn{1}{|c||}{order}   \\ \hline
3.125e$-$02 & 7.017e$-$06 & 4.69 & 3.420e$-$07 & 6.36 \\
1.563e$-$02 & 2.282e$-$07 & 4.94 & 2.646e$-$09 & 7.01 \\
7.813e$-$03 & 7.164e$-$09 & 4.99 & 2.123e$-$11 & 6.96 \\
3.906e$-$03 & 2.240e$-$10 & 5.00 & 1.706e$-$13 & 6.96 \\
\hline
\end{tabular}
\caption{Local error and deviation of the symmetrized defect-based error estimator for the fourth order integrator
from~\cite[\texttt{Emb 4/3 AK s}]{splithp} applied to~\eqref{eq:cnls}.\label{tab:conv2}}
\end{center}
\end{table}

\begin{table}[!ht]
\begin{center}
\begin{tabular}{||r||r|r||c|r||}
\hline
\multicolumn{1}{||c||}{$\tau$}    & \multicolumn{1}{c|}{global error} & \multicolumn{1}{|c||}{order}   &
\multicolumn{1}{c|}{error of corrected solution} &  \multicolumn{1}{|c||}{order}   \\ \hline
3.125e$-$02 & 7.894e$-$06 & 4.85 & 6.859e$-$07 & 5.97 \\
1.563e$-$02 & 4.035e$-$07 & 4.29 & 2.771e$-$09 & 7.95 \\
7.813e$-$03 & 2.471e$-$08 & 4.03 & 2.987e$-$11 & 6.54 \\
3.906e$-$03 & 1.537e$-$09 & 4.01 & 4.622e$-$13 & 6.01 \\
\hline
\end{tabular}
\caption{\rev{Global error and corrected solution for the fourth order integrator
from~\cite[\texttt{Emb 4/3 AK s}]{splithp} applied to~\eqref{eq:cnls}}.\label{tab:conv2b}}
\end{center}
\end{table}

\rev{\paragraph{Adaptive time-stepping}
The error estimators introduced in this paper are intended to be used as the basis for an adaptive
time-stepping procedure to enhance the efficiency. To illustrate this aspect, we show step-sizes
generated by the standard step-size selection strategy \cite{haireretal87}. We solve problem
(\ref{eq:cnls}) with the initial condition
$$ \psi(x,0)= \sum_{j=1}^2 \frac{a_j \ee^{-\ii b_j x}}{\cosh(a_j(x-c_j))}$$
with $a_1=a_2=2,\ b_1=1,\ b_2 =-3,\ c_1=5,\ c_2=-5,$ and a space discretization at
512 points on the interval $[-16,16].$  Time integration is effected by the integrator from \cite[\texttt{Emb 4/3 AK s}]{splithp}.
This example features two solitons which cross at $t\approx2.3$,
at which point the unsmooth solution demands smaller stepsizes.
If we prescribe a tolerance
of $10^{-10}$ on the local error, we obtain the stepsizes shown in Figure~\ref{fig:steps}.
It is found that the stepsizes indeed decrease in the region where the solitons cross,
which corresponds with the behavior observed for adaptive time-stepping based on
standard error estimators in \cite{koarl2paper}.
\begin{figure}[ht!]
\begin{center}
\includegraphics[width=0.8\textwidth]{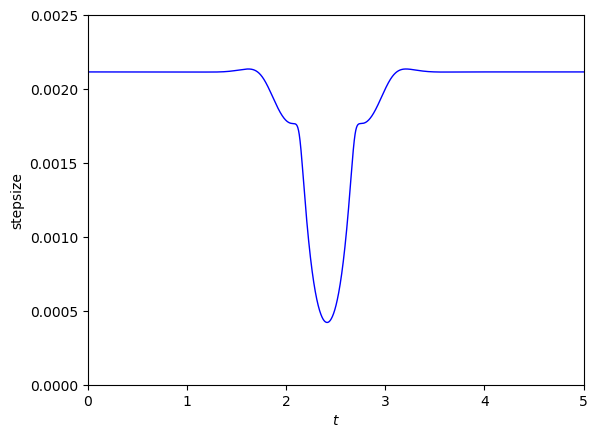}
\caption{\rev{Step-sizes generated by an adaptive strategy based on the symmetric error estimator
for the integrator from \cite[\texttt{Emb 4/3 AK s}]{splithp} for the problem (\ref{eq:cnls})
with crossing solitons.}\label{fig:steps}}
\end{center}
\end{figure}
}

\subsection{Rosen--Zener model}
As a second example, we solve a Rosen-Zener model from~\cite{blacastha17} by Magnus-type methods.
The associated Schr{\"o}dinger equation in the interaction
picture is given by
\begin{equation}\label{rosen}
\mathrm{i}\dot\psi(t) = H(t)\psi(t)
\end{equation}
with
\begin{equation} \label{zener}
\begin{aligned}
& H(t) = f_1(t) \sig_1 \otimes I_{k\times k} + f_2(t) \sig_2 \otimes R \in\mathbb{C}^{2k\times 2k},\quad k=50,\\
& \sig_1 = \left( \begin{array}{cc} 0 & 1 \\ 1 & 0 \end{array} \right),\qquad
  \sig_2 = \left( \begin{array}{rr} 0 & -\ii \\ \ii & 0 \end{array} \right),\\
& R = \mathrm{tridiag}(1,0,1) \in \mathbb{R}^{k\times k},\qquad f_1(t) = V_0 \cos(\omega t) \left( \cosh(t/T_0) \right)^{-1},\\
& f_2(t) = V_0 \sin(\omega t) \left( \cosh(t/T_0) \right)^{-1},
\qquad \omega=\tfrac12,\ T_0 =1,\ V_0 = 1,
\end{aligned}
\end{equation}
subject to the initial condition $\psi(0)=(1,\dots,1)^T.$

In Tables~\ref{tab:conv1rosen}--\ref{tab:conv2rosene}, we give the local errors and
deviation of the symmetrized error estimators for the test problem~\eqref{rosen}.
Table~\ref{tab:conv1rosen} gives the results for the exponential midpoint rule,
where the symmetrized defect can be evaluated exactly. Tables~\ref{tab:conv2rosenb}
and~\ref{tab:conv2rosen} give the empirical convergence orders for the commutator-free fourth order
Magnus-type integrator~\cite[\texttt{CF4:2} in Table~2]{alvfeh11} in conjunction with the symmetrized
defect-based error estimator,
evaluated by means of the Taylor variant
in Table~\ref{tab:conv2rosenb} and the Hermite variant
in Table~\ref{tab:conv2rosen}, respectively
(see~Figure~\ref{fig:algscfm}).
Finally, Table~\ref{tab:conv3rosen} gives the result
for the classical fourth order Magnus integrator, where the error estimator
is evaluated by means of the Hermite variant.
\rev{Tables~\ref{tab:conv1rosenb}, \ref{tab:conv2rosenc}, \ref{tab:conv2rosend} and
\ref{tab:conv2rosene} give the corresponding global errors on the interval $[0,1]$ of the basic solution
and of the solution corrected by the symmetric error estimate.}
In all cases, the theoretical results are
well reflected in the numerical experiments.

\begin{table}[!ht]
\begin{center}
\begin{tabular}{||r||r|r||c|r||}
\hline
\multicolumn{1}{||c||}{$\tau$}    & \multicolumn{1}{c|}{$\|\nL(\tau,u_0)\|_2$} & \multicolumn{1}{|c||}{order}   &
\multicolumn{1}{c|}{$\|{\nLLtildes}(\tau,u_0) - \nL(\tau,u_0)\|_2$} &  \multicolumn{1}{|c||}{order}   \\ \hline
1.250e$-$01 & 3.343e$-$03 & 2.97 & 7.157e$-$06 & 4.96\\
6.250e$-$02 & 4.198e$-$04 & 2.99 & 2.251e$-$07 & 4.99\\
3.125e$-$02 & 5.254e$-$05 & 3.00 & 7.047e$-$09 & 5.00\\
1.563e$-$02 & 6.569e$-$06 & 3.00 & 2.203e$-$10 & 5.00\\
7.813e$-$03 & 8.212e$-$07 & 3.00 & 6.885e$-$12 & 5.00\\
3.906e$-$03 & 1.026e$-$07 & 3.00 & 2.157e$-$13 & 5.00\\
\hline
\end{tabular}
\caption{Local error and deviation of the symmetrized defect-based error estimator for the second order exponential midpoint
rule applied to~\eqref{rosen}.\label{tab:conv1rosen}}
\end{center}
\end{table}

\begin{table}[!ht]
\begin{center}
\begin{tabular}{||r||r|r||c|r||}
\hline
\multicolumn{1}{||c||}{$\tau$}    & \multicolumn{1}{c|}{global error} & \multicolumn{1}{|c||}{order}   &
\multicolumn{1}{c|}{error of corrected solution} &  \multicolumn{1}{|c||}{order}   \\ \hline
5.000e$-$01 & 2.713e-01 &      & 7.652e-03 & \\
2.500e$-$01 & 6.618e-02 & 2.04 & 4.638e-04 & 4.04\\
1.250e$-$01 & 1.645e-02 & 2.01 & 2.880e-05 & 4.01\\
6.250e$-$02 & 4.106e-03 & 2.00 & 1.797e-06 & 4.00\\
3.125e$-$02 & 1.026e-03 & 2.00 & 1.123e-07 & 4.00\\
1.563e$-$02 & 2.565e-04 & 2.00 & 7.018e-09 & 4.00\\
\hline
\end{tabular}
\caption{\rev{Global error and corrected solution for the exponential midpoint rule applied to~\eqref{rosen}}.\label{tab:conv1rosenb}}
\end{center}
\end{table}

\begin{table}[!ht]
\begin{center}
\begin{tabular}{||r||r|r||c|r||}
\hline
\multicolumn{1}{||c||}{$\tau$}    & \multicolumn{1}{c|}{$\|\nL(\tau,u_0)\|_2$} & \multicolumn{1}{|c||}{order}   &
\multicolumn{1}{c|}{$\|{\nLLtildes}(\tau,u_0) - \nL(\tau,u_0)\|_2$} &  \multicolumn{1}{|c||}{order}   \\ \hline
\hline
5.000e$-$01 & 1.884e$-$03 & 4.78 & 5.854e$-$05 & 6.61 \\
2.500e$-$01 & 6.029e$-$05 & 4.97 & 4.875e$-$07 & 6.91 \\
1.250e$-$01 & 1.892e$-$06 & 4.99 & 3.868e$-$09 & 6.98 \\
6.250e$-$02 & 5.918e$-$08 & 5.00 & 3.033e$-$11 & 6.99 \\
3.125e$-$02 & 1.850e$-$09 & 5.00 & 2.373e$-$13 & 7.00 \\
\hline
\end{tabular}
\caption{Local error and deviation of the symmetrized defect-based error estimator for the fourth order CFM integrator
\cite[\texttt{CF4:2} in Table~2]{alvfeh11} applied to~\eqref{rosen},
defect evaluation by Taylor variant.
\label{tab:conv2rosenb}}
\end{center}
\end{table}

\begin{table}[!ht]
\begin{center}
\begin{tabular}{||r||r|r||c|r||}
\hline
\multicolumn{1}{||c||}{$\tau$}    & \multicolumn{1}{c|}{global error} & \multicolumn{1}{|c||}{order}   &
\multicolumn{1}{c|}{error of corrected solution} &  \multicolumn{1}{|c||}{order}   \\ \hline
5.000e$-$01 & 2.098e-03 &      & 5.330e-05 & \\
2.500e$-$01 & 1.212e-04 & 4.11 & 7.419e-07 & 6.17\\
1.250e$-$01 & 7.443e-06 & 4.03 & 1.126e-08 & 6.04\\
6.250e$-$02 & 4.632e-07 & 4.01 & 1.745e-10 & 6.01\\
3.125e$-$02 & 2.892e-08 & 4.00 & 2.768e-12 & 5.98\\
1.563e$-$02 & 1.807e-09 & 4.00 & 1.175e-13 & 4.56\\
\hline
\end{tabular}
\caption{\rev{Global error and corrected solution for the fourth order CFM integrator
\cite[\texttt{CF4:2} in Table~2]{alvfeh11} applied to~\eqref{rosen},
defect evaluation by Taylor variant.\label{tab:conv2rosenc}}}
\end{center}
\end{table}

\begin{table}[!ht]
\begin{center}
\begin{tabular}{||r||r|r||c|r||}
\hline
\multicolumn{1}{||c||}{$\tau$}    & \multicolumn{1}{c|}{$\|\nL(\tau,u_0)\|_2$} & \multicolumn{1}{|c||}{order}   &
\multicolumn{1}{c|}{$\|{\nLLtildes}(\tau,u_0) - \nL(\tau,u_0)\|_2$} &  \multicolumn{1}{|c||}{order}   \\ \hline
5.000e$-$01 & 1.884e$-$03 & 4.78 & 4.008e$-$05 & 6.64 \\
2.500e$-$01 & 6.029e$-$05 & 4.97 & 3.277e$-$07 & 6.93 \\
1.250e$-$01 & 1.892e$-$06 & 4.99 & 2.584e$-$09 & 6.99 \\
6.250e$-$02 & 5.918e$-$08 & 5.00 & 2.023e$-$11 & 7.00 \\
3.125e$-$02 & 1.850e$-$09 & 5.00 & 1.583e$-$13 & 7.00 \\
\hline
\end{tabular}
\caption{Local error and deviation of the symmetrized defect-based error estimator for the fourth order CFM integrator
\cite[\texttt{CF4:2} in Table~2]{alvfeh11} applied to~\eqref{rosen},
defect evaluation by Hermite variant.
\label{tab:conv2rosen}}
\end{center}
\end{table}

\begin{table}[!ht]
\begin{center}
\begin{tabular}{||r||r|r||c|r||}
\hline
\multicolumn{1}{||c||}{$\tau$}    & \multicolumn{1}{c|}{global error} & \multicolumn{1}{|c||}{order}   &
\multicolumn{1}{c|}{error of corrected solution} &  \multicolumn{1}{|c||}{order}   \\ \hline
5.000e$-$01 & 2.098e-03 &      & 3.203e-05 & \\
2.500e$-$01 & 1.212e-04 & 4.11 & 4.402e-07 & 6.19\\
1.250e$-$01 & 7.443e-06 & 4.03 & 6.702e-09 & 6.04\\
6.250e$-$02 & 4.632e-07 & 4.01 & 1.041e-10 & 6.01\\
3.125e$-$02 & 2.892e-08 & 4.00 & 1.676e-12 & 5.96\\
1.563e$-$02 & 1.807e-09 & 4.00 & 1.052e-13 & 3.99\\
\hline
\end{tabular}
\caption{\rev{Global error and corrected solution for the fourth order CFM integrator
\cite[\texttt{CF4:2} in Table~2]{alvfeh11} applied to~\eqref{rosen},
defect evaluation by Hermite variant.\label{tab:conv2rosend}}}
\end{center}
\end{table}

\begin{table}[!ht]
\begin{center}
\begin{tabular}{||r||r|r||c|r||}
\hline
\multicolumn{1}{||c||}{$\tau$}    & \multicolumn{1}{c|}{$\|\nL(\tau,u_0)\|_2$} & \multicolumn{1}{|c||}{order}   &
\multicolumn{1}{c|}{$\|{\nLLtildes}(\tau,u_0) - \nL(\tau,u_0)\|_2$} &  \multicolumn{1}{|c||}{order}   \\ \hline
5.000e$-$01 & 4.788e$-$03 & 4.56 & 1.214e$-$04 & 6.13 \\
2.500e$-$01 & 1.618e$-$04 & 4.89 & 1.126e$-$06 & 6.75 \\
1.250e$-$01 & 5.154e$-$06 & 4.97 & 9.201e$-$09 & 6.94 \\
6.250e$-$02 & 1.618e$-$07 & 4.99 & 7.269e$-$11 & 6.98 \\
3.125e$-$02 & 5.064e$-$09 & 5.00 & 5.693e$-$13 & 7.00 \\
\hline
\end{tabular}
\caption{Local error and deviation of the symmetrized defect-based error estimator
for the fourth order classical Magnus integrator~\eqref{MC-4}
applied to~\eqref{rosen}, defect evaluation by Hermite variant.
\label{tab:conv3rosen}}
\end{center}
\end{table}

\begin{table}[!ht]
\begin{center}
\begin{tabular}{||r||r|r||c|r||}
\hline
\multicolumn{1}{||c||}{$\tau$}    & \multicolumn{1}{c|}{global error} & \multicolumn{1}{|c||}{order}   &
\multicolumn{1}{c|}{error of corrected solution} &  \multicolumn{1}{|c||}{order}   \\ \hline
5.000e-01 & 6.957e-03 &      & 1.536e-04 & \\
2.500e-01 & 4.362e-04 & 4.00 & 2.452e-06 & 5.97\\
1.250e-01 & 2.728e-05 & 4.00 & 3.853e-08 & 5.99\\
6.250e-02 & 1.705e-06 & 4.00 & 6.029e-10 & 6.00\\
3.125e-02 & 1.066e-07 & 4.00 & 9.419e-12 & 6.00\\
1.563e-02 & 6.662e-09 & 4.00 & 1.688e-13 & 5.80\\\hline
\end{tabular}
\caption{\rev{Global error and corrected solution for the
fourth order classical Magnus integrator~\eqref{MC-4} applied to~\eqref{rosen},
defect evaluation by Hermite variant.\label{tab:conv2rosene}}}
\end{center}
\end{table}

\section{Conclusion} \label{sec:concl}

We have discussed a symmetrized defect-based estimator for self-adjoint
time discretizations of nonlinear evolution equations. We have introduced
the general construction principle extending the ideas from~\cite{auzingeretal18a},
and have elaborated the algorithms for an implicit Runge-Kutta method,
for splitting methods and for exponential Magnus-type integrators
for time-dependent linear problems. We have proven
that the deviation of the estimated error from the true error is two
orders in the step-size smaller than the basic integrator, and illustrated
the theoretical result for two examples solved by either splitting methods
or exponential Magnus-type integrators of different orders.

It can be expected that in adaptive simulations, where choice of the step-size
is delicate, the improved accuracy of the error estimator may add to
the reliability and efficiency of the integrator.
However, this topic exceeds
the scope of the present work and will be explored elsewhere.
\rev{Here, we have confined ourselves to a numerical illustration that
our error estimators induce adaptive step-sizes commensurate with the
solution behavior.}
Note, moreover, that the numerical approximation
based on a scheme of order $ p $ and corrected
by our error estimator (see~\eqref{Sdach-nonlin})
is very close to self-adjoint
and has improved convergence order $ p+2 $
(see Theorem~\ref{Sdach-results}),
thus providing a nearly self-adjoint higher order
approximation at moderate computational cost.
Since the additive correction is of high order,
no stability problems will arise for
the corrected scheme~\eqref{Sdach-nonlin}.

\section*{Acknowledgements}
This work was supported in part by the Vienna Science and Technology Fund (WWTF) [grant number MA14-002]
and the Austrian Science Fund (FWF) [grant number P 30819-N32]. \rev{We thank D. Haberlik, student at
TU Wien, for contributing some of the numerical results, and}
M.~Brunner, student at TU Wien, for contributing Figure~\ref{fig:GALS}.


\end{document}